\documentclass[12pt]{article}
\usepackage[utf8]{inputenc}
\usepackage[T1]{fontenc}
\usepackage{lmodern}
\usepackage[margin=2cm]{geometry}
\usepackage[english]{babel}
\usepackage{amsmath, amssymb, amsthm}
\usepackage{hyperref}
\usepackage[all]{xy}
\usepackage{tikz}
\usepackage{xypic}

\newtheorem{thm}{Theorem}[section]
\newcounter{thmL}
\newtheorem{cor}[thm]{Corollary}

\newtheorem{defi}[thm]{Definition}

\newtheorem{lem}[thm]{Lemma}

\newtheorem{prop}[thm]{Proposition}

\newcommand*{\QED}{\hfill\ensuremath{\square}}
\newenvironment{demo}{\noindent \textit{Proof.}}{\hfill\QED\\}
\newenvironment{demoKW}{\noindent \textit{Proof of Theorem \ref{KW}.}}{\hfill\QED\\}

\newenvironment{rmq}{\noindent \textbf{Remark.}}{}

\renewcommand*{\Re}{\mathfrak{Re}}

\DeclareMathOperator{\Gal}{Gal}
\DeclareMathOperator{\Frob}{Frob}
\DeclareMathOperator{\ord}{ord}
\DeclareMathOperator{\Var}{Var}
\DeclareMathOperator{\Span}{Span}

\newcommand*{\C}{\mathbb C}
\newcommand*{\E}{\mathbb E}
\newcommand*{\F}{\mathbb F}
\renewcommand*{\H}{\mathbb H}

\newcommand*{\N}{\mathbb N}

\renewcommand*{\P}{\mathbb P}
\newcommand*{\Q}{\mathbb Q}
\newcommand*{\R}{\mathbb R}

\newcommand*{\T}{\mathbb T}

\newcommand*{\Z}{\mathbb Z}

\begin{document}

\title{\textbf{Explicit Kronecker-Weyl theorems and applications to prime number races}}
\author{Alexandre Bailleul}
\date{}

\maketitle

\abstract{We prove explicit versions of the Kronecker-Weyl theorems, both in a discrete and a continuous settings, without any linear independence hypothesis. As an application, we propose an alternative approach to problems concerning asymptotic densities in prime number races, over number fields and over function fields in one variable over finite fields, in the language of random variables. Our approach allows us to prove new results on the existence and positivity of some of those densities, which, in the case of races over function fields, do not require any linear independence hypothesis.}

\section{Introduction}

\subsection{Context}

The Kronecker-Weyl theorem is an important result of harmonic analysis, related to ergodic theory, and which has been used to study many arithmetical problems of statistical nature. It is both a multidimensional generalization of Weyl's famous result on the equidistribution of the fractional parts of $n\alpha$ ($n \in \N$), when $\alpha$ is an irrational number, and a generalization of Kronecker's density result on the torus. More specifically, let $\theta_1, \dots, \theta_n$ be real numbers, then the one-parameter subgroup $$\Gamma := \left\{\left(e^{i\theta_1x}, \dots, e^{i\theta_nx}\right) \mid x \in \R\right\}$$ is equidistributed in a subtorus inside the $n$-dimensional torus $$\T^n := \{(z_1, \dots, z_n) \in \C^n \mid \forall i \in \{1, \dots, n\}, |z_i|=1\}$$ with respect to its Haar measure $\mathrm{d}\mu$. In other words, the topological closure $\overline{\Gamma}$ of $\Gamma$ in $\T^n$ is a closed subgroup of $\T^n$ with Haar measure $\mathrm{d}\mu$ and for every continuous function $f : \T^n \rightarrow \C$, one has $$\frac{1}{X} \int_0^X f\left(e^{i\theta_1x}, \dots, e^{i\theta_nx}\right) \,\mathrm{d}x \underset{X \to +\infty}{\longrightarrow} \int_{\overline{\Gamma}} f \,\mathrm{d}\mu.$$ Most relevant to the present work is the additional information that $\overline{\Gamma}$ is a $m$-dimensional torus, where $m$ is the dimension of the $\Q$-vector space spanned by $\theta_1, \dots, \theta_n$. In particular, if $\theta_1, \dots, \theta_n$ are $\Q$-linearly independent, then $\overline{\Gamma} = \T^n$, so we obtain Kronecker's density result in a strong form (in the sense that equidistribution holds), and when $n=1$, this is exactly Weyl's equidistribution result.

There exists a discrete version of the Kronecker-Weyl Theorem, in which we consider the discretely-parametrized subgroup $$\Gamma := \left\{\left(e^{i\theta_1X}, \dots, e^{i\theta_nX}\right) \mid X \in \Z\right\}.$$ In this case, integrals are replaced by sums and we require the real numbers $\theta_1, \dots, \theta_n$ to be $\Q$-linearly independent with $\pi$. The reason for this is clear, since $e^{iq\pi X}$ assumes discrete values in $\T$ when $X$ ranges over the integers, and $q$ is a rational number. Usually, both the continuous and the discrete versions of the Kronecker-Weyl Theorem are proved using abstract harmonic analysis (see \cite[Theorem 2.2.5]{Que} or \cite[Theorem 4.2]{Dev} for instance). In this paper, we give an elementary proof of a general version of Kronecker-Weyl's result, both in the discrete and the continuous case, in which we explicitly construct the set in which $\Gamma$ equidistributes. We insist on the fact that no hypothesis of linear independance is required in our result (see Corollary Theorem \ref{Prob} and Corollary \ref{CProb}). We note that an elementary and explicit proof of the continuous version of the Kronecker-Weyl theorem was given in the arXiv version \cite{MNa} of \cite{MN}, Appendix B.\\

The Kronecker-Weyl theorem is at the heart of the modern approach (initiated by Rubinstein and Sarnak in \cite{RuSa}) to the study of so-called "prime number races", which consists in investigating the properties of the set $$\mathcal{P}_{q; a_1, \dots, a_D} := \{x \geq 2 \mid \pi(x;q,a_1) > \pi(x;q,a_2) > \dots > \pi(x;q,a_D)\},$$ where $\pi(x; k ,c)$ is the number of primes $p \leq x$ such that $p \equiv c$ mod $k$. Here, the invertible classes $a_1, \dots, a_D$ are called the \textit{contestants} of the prime number race. Assuming the Generalized Riemann Hypothesis, Rubinstein and Sarnak proved that functions of the form $$E_{q; a_1, \dots, a_D} : y \mapsto \left(\pi(e^y;q,a_1) - \frac{\text{Li}(e^y)}{\varphi(q)}, \dots, \pi(e^y;q,a_D) - \frac{\text{Li}(e^y)}{\varphi(q)}\right)$$ admit limiting distributions, according to the following definition.

\begin{defi} Let $E : \R^+ \longrightarrow \R^D$. We say $E$ admits a limiting distribution $\mu$ when $\mu$ is a Borel probability measure on $\R^d$ such that for any bounded continuous function $f : \R^D \longrightarrow \R$, one has $$\frac{1}{X} \int_0^X f(E(y)) \,\mathrm{d}y \underset{X \to +\infty}{\longrightarrow} \int_{\R^D} f \,\mathrm{d}\mu.$$
\end{defi}

From there, and assuming the linear independence over $\Q$ of the non-negative imaginary parts of non-trivial zeros of Dirichlet $L$-functions mod $q$ (an hypothesis called the Grand Simplicity Hypothesis or GSH in \cite{RuSa}), Rubinstein and Sarnak proved that sets of the form $$\{x \geq 2 \mid \pi(e^x;q,a_1) > \pi(e^x;q,a_2) > \dots > \pi(e^x;q,a_D)\}$$ admit natural densities strictly between $0$ and $1$.

The discrete version of the Kronecker-Weyl theorem has been used initially by Cha \cite{Cha}, followed by other authors (\cite{ChaIm}, \cite{CFJ}, \cite{DevMen}), to study other kinds of prime number races over function fields, still assuming some form of linear independence between the zeros of the corresponding $L$-functions.

Further works aimed at weakening those linear independence hypotheses. In \cite{MN}, Martin and Ng introduced the notions of \textit{exhaustivity}, \textit{weak inclusiveness} and \textit{inclusiveness} of prime number races, focusing on the unboundedness (resp. the existence of the logarithmic densities, resp. the positivity of the logarithmic densities) of sets of the form $\mathcal{P}_{q; a_1, \dots, a_D}$. They introduced the notion of \emph{self-sufficient zero} (\cite[Definition 1.3]{MN}), and assuming various hypotheses on the existence of such zeros, proved the weak inclusivenes or inclusiveness of the corresponding prime number races. These assumptions were weakened by Devin in \cite{Dev}, who even studied the regularity of the corresponding limiting distributions in a more general setting.\\

Among the aforementioned properties we will mostly be interested in
weak inclusiveness and inclusiveness. We recall the definition of these notions in the context of general prime number races.

\begin{defi} \label{Incl} Let $a_1, \dots, a_D$ be contestants in a prime number race. We say the race between $a_1, \dots, a_D$ is weakly inclusive if for every permutation $\sigma$ of $\{1, \dots, D\}$, the set $\{X \geq 2 \mid \Pi(X, a_{\sigma(1)}) > \dots > \Pi(X, a_{\sigma(D)})\}$ admits a natural density, where the $\Pi( \cdot, a)$ are the corresponding (rescaled) prime counting functions. We say the prime number race is inclusive if moreover those densities are positive.
\end{defi}

The prime number races referred to in the above definition will be of two types in the present work. First, a prime number race over a number field denotes the race between unramified prime ideals in a Galois extension $L/K$ of number fields with given Frobenius automorphisms in $\Gal(L/K)$, as was first suggested in \cite[Section 5]{RuSa} and first studied in \cite[Chapter 5]{Ng}. In that case, the contestants are distinct conjugacy classes $C_1, \dots, C_D$ of $\Gal(L/K)$, and the rescaled prime counting functions are \begin{align*}\Pi(X, C_i) &:= \frac{\pi(e^X; L/K, C_i)}{\#C_i}\\ &= \frac{1}{\#C_i} \#\{\mathfrak{p} \text{ prime ideal of } K \text{ unramified in } L \mid N\mathfrak{p} \leq e^X, \text{Frob}_{\mathfrak{p}} = C_i\}.\end{align*} That the variable has to be changed to $e^X$ comes from the shape of the explicit formulas involved (see \cite[Chapter 5]{Ng}).

Second, a prime number race over a function field denotes the race between unramified prime divisors in a Galois extension $L/K$ of functions fields in one variable over a finite field with given Frobenius automorphisms in $\Gal(L/K)$. In that case, the contestants are distinct conjugacy classes $C_1, \dots, C_D$ of $\Gal(L/K)$, and the rescaled prime counting functions are \begin{align*}\Pi(X, C_i) &:= \frac{\pi(X; L/K, C_i)}{\#C_i}\\ &= \frac{1}{\#C_i} \#\{P \text{ prime divisor of } K \text{ unramified in } L \mid \deg P = X, \text{Frob}_{P} = C_i\}.\end{align*} One can also consider functions counting prime divisors of $K$ with a given Frobenius automorphism of degree less than $X$, instead of equal to $X$, as was studied in \cite{ChaIm}. Of course we recover the classical case of Rubinstein and Sarnak for primes in arithmetic progressions \cite{RuSa} in the number field case by considering an appropriate cyclotomic extension of $\Q$. Similarly we recover the races between irreducible polynomials in arithmetic progressions of \cite{Cha} by considering an appropriate Carlitz extension of $\F_q(T)$.

More general races have been studied in the literature, for instance the race between prime quadratic residues and prime non-quadratic residues modulo an integer $q$ in \cite{RuSa}, the race between $\pi(x)$ and $\text{Li}(x)$ in \cite{RuSa} and \cite{ANS}, the race between products of $k$ irreducible polynomials over a finite field in \cite{DevMen}, the race between the number of points on the reduction modulo good primes of elliptic curves in \cite{CFJ} and many more. In any case, it is clear to which category each of those races should belong, either over number fields or function fields. Our general results can be applied to those situations as well.\\

The first step in studying a prime number race over a function field is to write an explicit formula, \textit{i.e.} express the corresponding prime counting functions as sums involving $e^{i \theta_1 X}, \dots, e^{i \theta_r X}$, where $\theta_1, \dots, \theta_r$ are the positive arguments (between $0$ and $\pi$) of the inverse zeros of the corresponding rational $L$-functions. As an application of our version of the discrete Kronecker-Weyl theorem, we give sufficient conditions for the existence and for the positivity of the natural densities relevant to those kinds of races. Recently, Devin (\cite{Dev2}) studied the question of the existence of those densities, and provided sufficient conditions on the coefficients of the functions involved in the explicit formulas. Our approach is transverse to hers, as we give conditions on the functions themselves. Our approach is also considerably more elementary, as Devin relies on multiple tools of harmonic analysis to deduce that "ties have density zero" in such races. We avoid the use of such techniques thanks to our approach based on random variables and our key Lemma \ref{Atom}.

In the case of prime number races over number fields, the situation is technically more complicated, since explicit formulas for the (rescaled) prime counting functions involve infinite series in $e^{i \theta_1 t}, e^{i \theta_2 t}, \dots$ where $\theta_1, \theta_2, \dots$ are the positive imaginary parts of non-trivial zeros of the corresponding $L$-functions. Functions of this shape are often called almost periodic functions. There are different classes of almost periodic functions, depending on the way they can be approximated by trigonometric polynomials. That the remainder in the prime number theorem is almost-periodic dates back to at least Wintner (\cite{Win}, see also more recent works \cite{ANS} and \cite{KR} for more general remainders). The class of functions which is most relevant to us is the (large) class of Besicovitch almost periodic functions, called $B^1$-almost periodic functions. The existence of the limiting distributions of such functions is shown for example in \cite[Theorem 2.9]{ANS}. See also \cite[Theorem 4.1]{Ble} for a similar proof in a slightly larger space than $B^1$. As an application of our version of the Kronecker-Weyl theorem, we give a more precise description of this limiting distribution, under various hypotheses on the almost periods $\theta_n$ of the $B^1$-almost periodic function which is being studied (Corollary \ref{L2}). We also give a new proof of a recent result of Devin, giving sufficient conditions for the existence of the densities associated to $B^1$-almost periodic functions. Our approach is again more elementary and does not require the use of abstract harmonic analysis (see Corollary \ref{Probinf}), though we appeal to Lévy's criterion for weak convergence of measures. We give an application to the existence of the densities involved in prime number races over number fields in Theorem \ref{NField}.

\subsection{Organization of the paper and statements of results}

The paper is organized as follows. In the first part, we prove explicit versions of the discrete and the continuous Kronecker-Weyl theorems for real numbers $\theta_1, \dots, \theta_r$.

\begin{thm} Let $\Gamma_{\theta} = \left\{\left(e^{i\theta_1X}, \dots, e^{i\theta_rX}\right) \mid X \in \Z\right\}$. Then $\Gamma_{\theta}$ is equidistributed in $\overline{\Gamma_{\theta}} = \bigcup_{a=0}^{d-1} \nu_{\theta}^aH_{\theta} = \langle\nu_{\theta}\rangle H_{\theta}$, where $\nu_{\theta} = \left(e^{i \theta_1}, \dots, e^{i \theta_r}\right)$, and $$H_{\theta} = \left\{\left(z_1^d, \dots, z_m^d, \prod_{k=1}^m z_k^{h_{k, m+1}}, \dots, \prod_{k=1}^m z_k^{h_{k, r}}\right) \mid (z_1, \dots, z_m) \in \T^m\right\} \subset \T^r,$$ that is for every continuous $f : \T^r \longrightarrow \C$ one has $$\frac{1}{X} \sum_{n \leq X} f\left(e^{i\theta_1 n}, \dots, e^{i \theta_r n}\right) \underset{X \to +\infty}{\longrightarrow} \int_{\overline{\Gamma_{\theta}}} f \,\mathrm{d}\mu_{\theta}.$$ The measure $\mu_{\theta}$ with respect to which $\Gamma_{\theta}$ is equidistributed is $\frac{1}{d}\sum_{a=0}^{d-1} \mu_a$, where $\mu_a$ is the pushforward of the normalized Haar measure $\mu_{H_{\theta}}$ of $H_{\theta}$ to $\nu_{\theta}^aH_{\theta}$.
\end{thm}

The definitions of $m, d$ and $h_{k, j}$ are given at the beginning of Section \ref{Disc}.

\begin{thm} The one-parameter subgroup $$\Gamma_{\theta} = \left\{\left(e^{i\theta_1 y}, \dots, e^{i\theta_r y}\right) \mid y \in \R\right\}$$ is equidistributed in $$H_{\theta} = \left\{\left(z_1^d, \dots, z_m^d, \dots, \prod_{k=1}^m z_k^{h_{k,j}}, \dots\right) \mid (z_1, \dots, z_m) \in \T^m\right\},$$ that is for every continuous $f : \T^r \longrightarrow \C$ one has $$\frac{1}{X} \int_0^X f\left(e^{i\theta_1 y}, \dots, e^{i \theta_r y}\right) \,\mathrm{d}y \underset{X \to +\infty}{\longrightarrow} \int_{H_{\theta}} f \,\mathrm{d}\mu_{H_{\theta}}$$ where $\mu_{H_{\theta}}$ is the normalized Haar measure on $H_{\theta}$.
\end{thm}

For this statement, the definitions of $m$, $d$ and $h_{k, j}$ are given at the beginning of Section \ref{Cont}.\\


Those theorems are our starting point to obtain the existence of asymptotic densities for sets defined by strict inequalities between certain types of functions without any linear independence hypothesis. In doing so, we prove Lemma \ref{Atom} which allows us to bypass technical results from harmonic analysis to prove that the limiting distributions we consider do not admit atoms in the non-degenerate discrete case.

\begin{thm}
For $1 \leq j \leq D$, let $f_j \in \C(X_1, \dots, X_r)$ be real-valued and without pole on $\T^r$ and let $F_j : t \mapsto f_j\left(e^{i\theta_1 t}, \dots, e^{i\theta_rt}\right)$. Then we have $$\frac{1}{X} \#\left\{n \leq X \mid F_1(n) > \dots > F_D(n)\right\} \underset{X \to +\infty}{\longrightarrow} \frac{1}{d} \sum_{a=0}^{d-1} \P(f_1(\nu_{\theta}^aZ_{\theta}) > \dots > f_D(\nu_{\theta}^aZ_{\theta}))$$ for some explicit random variable $Z_{\theta}$ on $\T^r$.
\end{thm}

The definition of $Z_{\theta}$ in this discrete case is given in (\ref{VA}).

\begin{thm} Let $f_1, \dots, f_D \in \C(X_1, \dots, X_r)$ be real-valued and without poles on $\T^r$ and let $F_j : t \mapsto f_j\left(e^{i\theta_1 t}, \dots, e^{i\theta_rt}\right)$. Then we have $$\frac{1}{X} \int_0^X \mathbf{1}_{x_1 > \dots > x_D}\left(F_1(y), \dots, F_D(y)\right) \,\mathrm{d}y \underset{X \to +\infty}{\longrightarrow} \P\left(f_1(Z_{\theta}) > \dots > f_D(Z_{\theta})\right)$$ for some explicit random variable $Z_{\theta}$ on $\T^r$.
\end{thm}

The definition of $Z_{\theta}$ in this continuous case is given in (\ref{VAc}).\\

We then generalize our methods in the continuous case, in the presence of infinitely many real numbers $\theta_1, \theta_2, \dots$. In this context, we prove that $B^1$-almost periodic functions admit limiting distributions (Theorem \ref{Helly}). We then tackle the problem of the existence of asymptotic densities for sets defined by strict inequalities between $B^1$-almost periodic functions, assuming suitable hypotheses (Proposition \ref{Infineg} and Theorem \ref{Probinf}).\\ 

The second part of the paper is devoted to applications. First, we apply our results in Section \ref{Erreur} to study the asymptotic densities associated to functions with an extra error term (Theorems \ref{Reste} and \ref{ResteC}), as those are the kind of functions appearing in explicit formulas for prime number races. For example, using the explicit formula \cite[(5.12)]{Ng} and our Theorem \ref{ResteC}, we get the following result for prime number races over number fields.

\begin{thm} \label{NField} Let $L/K$ be a Galois extension of number fields, with Galois group $G$. Assume $\zeta_L$ satisfies the Riemann Hypothesis. Let $\Theta$ be the set of positive imaginary parts of the non-trivial zeros of Artin $L$-functions attached to irreducible complex characters of $G$, and assume that $\Span_{\Q} \Theta = \Span_{\Q} T \oplus \Span_{\Q}(\Theta \setminus T)$ for some non-empty finite subset $T$ of $\Theta$. Let $C_1, \dots, C_D$ be distinct conjugacy classes of $G$. If for every $1 \leq j \leq D-1$, there exists $\theta \in T$ such that $$\sum_{\chi \neq \chi_0} \ord_{s=\frac{1}{2} + i \theta} L(s, \chi) \left(\chi(C_j) - \chi(C_{j+1})\right) \neq 0$$ then the logarithmic density $$\delta(L/K; C_1, \dots, C_D) := \lim_{X \to +\infty} \frac{1}{X} \int_2^X \mathbf{1}_{\pi_{C_1}(e^t) > \dots > \pi_{C_D}(e^t)} \,\mathrm{d}t$$ exists.
\end{thm}

\begin{rmq} \begin{itemize}
\item[i)] Using the unconditional explicit formula from \cite[Corollary 3.10]{FiJo}, we could provide a similar statement for the existence of the above logarithmic density, under a suitable hypothesis of non-vanishing coefficient as above and without assuming the Riemann Hypothesis for $\zeta_L$.

\item[ii)] Because of the special properties of Artin $L$-functions with respect to induction of characters, one could state a linear independence hypothesis about the set of zeros of Artin $L$-functions attached to irreducible complex characters of $G^+$ instead, where $G^+$ is the Galois group of the Galois closure of $L$ over $\Q$. That this is a more natural set of zeros to consider was noted by the author, and used for the first time in \cite{FiJo}.
\end{itemize}
\end{rmq}

In Section \ref{Noncrit}, we give general criteria for the positivity of those asymptotic densities (Propositions \ref{PosDisc} and \ref{PosCont}), with inclusiveness of prime number races over function fields in mind. In particular, as long as the main terms in explicit formulas for the prime counting functions satisfy certain inequalities for at least one value of $n$, then the corresponding prime counting functions do so for a positive proportion of $n$ (a $\liminf$ in general).\\

Then in Section \ref{Mom1} we study the first two moments of the limiting distributions appearing in the discrete case, as those can be used, for example with Chebyshev's inequality, to study the asymptotic density considered in prime number races over function fields. In Section \ref{Momente}, we do the same in the infinite-dimensional continuous case and we also give a description of the limiting distribution under a weak linear independence assumption (Theorem \ref{L2}).\\

Finally, in Section \ref{2} we apply those general results to the concrete problem of studying prime divisor races in geometric Galois extensions of function fields (in one variable) over finite fields. In particular, we are able to study an example of prime divisor race in which the usual linear independence hypothesis fails to hold.\\

It seems a reference to a proof of the discrete version of the Kronecker-Weyl theorem is hard to find in a published form so we provide a proof in an appendix. We borrowed the proof to P. Humphries' Masters thesis \cite{Hum}. For a proof of the general continuous version, see \cite[Theorem 4.2]{Dev}.\\

\textbf{Notations.} Some notations are introduced at various places and then used many times throughout the text, we gather them here for reference. In the discrete case, the element $\nu_{\theta}$ is defined in (\ref{nu}) and denotes $\left(e^{i\theta_1}, \dots, e^{i \theta_r}\right)$. The quantities $d$ and $h_{k, j}$ coming from linear dependence relations are introduced in (\ref{d}) and (\ref{h}) in Section \ref{Disc} and (\ref{dc}) and (\ref{hc}) in Section \ref{Cont}, depending on the context. The subgroup $H_{\theta} \subset \T^r$ arising in the non-degenerate case is defined in (\ref{H}) in the discrete case and in (\ref{Hc}) in the continuous case.  The random variables $Z_{\theta}$ are introduced, depending on the context, in formula (\ref{VA}) (discrete case) or formula (\ref{VAc}) (continuous case).

\section{Explicit Kronecker-Weyl theorems}

\label{1}

We begin by recalling the definition of equidistribution that we are going to use throughout the text.

\begin{defi} Let $(z_n)_{n \in \N}$ be a sequence of elements of $\T^r$ and $H$ a closed subgroup of $\T^r$, and let $\mu_H$ be the Haar measure of $H$. We say $(z_n)_{n \in \N}$ is equidistributed in $H$ with respect to the measure $\mu_H$ if for every continuous $f : \T^r \rightarrow \C$, one has $$\frac{1}{X} \sum_{n \leq X} f(z_n) \underset{X \to +\infty}{\longrightarrow} \int_{H} f \,\mathrm{d}\mu_H.$$
\end{defi}

This definition is equivalent to the weak convergence of the measures $\frac{1}{X}\sum_{n \leq X} \delta_{z_n}$ to the measure $\mu_H$, where $\delta_{z_n}$ is the Dirac measure at $z_n$. In what follows, we identify such a sequence with the set $Z := \{z_n \mid n \in \N\}$. This is a slight abuse of notation since the set $Z$ itself does not keep track of the numbering of the sequence. Note that if $Z$ is equidistributed in $H$, then $Z$ is in particular dense in $H$, so that $\overline{Z} = H$. The following is a weak version of the discrete Kronecker-Weyl theorem, which will be enough for our purpose of proving an explicit strong version.

\begin{thm}[Discrete Kronecker-Weyl theorem]\label{KW} Let $\theta_1, \dots, \theta_r$ be real numbers such that $\{\theta_1, \dots, \theta_r, \pi\}$ is linearly independent over $\Q$. Then the set $$\Gamma := \left\{\left(e^{i\theta_1X}, \dots, e^{i\theta_rX}\right) \mid X \in \Z\right\}$$ is equidistributed in $\T^r$ (with respect to its Haar measure).
\end{thm}

A proof is given in an appendix for reference. Notice the linear independence assumption with $\pi$. The goal of the following sections is to prove a precise version of this result (and of its continuous analog) with no assumption of linear independence (even with $\pi$) and with a description of the subset of $\T^r$ in which $\Gamma$ is equidistributed. We also do so by elementary means, while many proofs of the full Kronecker-Weyl theorem use abstract harmonic analysis (namely Pontryagin duality and Poisson summation formula).

\subsection{The discrete case}

\label{Disc}

In this section, consider real numbers $\theta_1, \dots, \theta_r$ and write $\theta = (\theta_1, \dots, \theta_r)$. Set $m+1 = \dim \Span_{\Q}(\pi, \theta_1, \dots, \theta_r)$. Up to reindexing the $\theta_i$'s, extract a basis $\{2\pi, \theta_1, \dots, \theta_m\}$ of the $\Q$-vector space $\Span_{\Q}(\pi, \theta_1, \dots, \theta_r)$ and write the decomposition of $\theta_{m+1}, \dots, \theta_r$ in this basis as \begin{equation}\label{coeff}
\theta_j = 2 \pi c_j + \sum_{k=1}^m b_{k,j} \theta_k \text{ for } m+1 \leq j \leq r
\end{equation} with $c_j, b_{k,j} \in \Q$. Note that if $m=0$ then all $\theta_i$'s are rational multiples of $\pi$ and the decomposition above reduces to $\theta_i = 2\pi c_i$ for $1 \leq i \leq r$. In that case we say we are in the \textit{degenerate case}, and otherwise in the \textit{non-degenerate case}. 

We let \begin{equation}\label{d}d := \mathrm{lcm}(\{\text{denominators of all } c_j \text{ and } b_{k,j}\}),\end{equation} so that \begin{equation}\label{h} l_j := dc_j \in \Z, h_{k,j} := db_{k,j} \in \Z.\end{equation} Finally, let \begin{equation}\label{nu}
\nu_{\theta} = \left(e^{i\theta_1}, \dots, e^{i\theta_r}\right).
\end{equation}

\begin{thm} \label{NKW} Let $\Gamma_{\theta} = \left\{\left(e^{i\theta_1X}, \dots, e^{i\theta_rX}\right) \mid X \in \Z\right\}$. Then $\Gamma_{\theta}$ is equidistributed in $\overline{\Gamma_{\theta}} = \bigcup_{a=0}^{d-1} \nu_{\theta}^aH_{\theta} = \langle\nu_{\theta}\rangle H_{\theta}$, and \begin{equation}\label{H} H_{\theta} = \left\{\left(z_1^d, \dots, z_m^d, \prod_{k=1}^m z_k^{h_{k, m+1}}, \dots, \prod_{k=1}^m z_k^{h_{k, r}}\right) \mid (z_1, \dots, z_m) \in \T^m\right\} \subset \T^r\end{equation} with Haar measure $\mu_{H_{\theta}}$. The measure $\mu_{\theta}$ with respect to which $\Gamma_{\theta}$ is equidistributed is $\frac{1}{d}\sum_{a=0}^{d-1} \mu_a$, where $\mu_a$ is the pushforward of $\mu_{H_{\theta}}$ to $\nu_{\theta}^aH_{\theta}$.
\end{thm}

Note that when $m=0$, $H_{\theta}$ reduces to the trivial subgroup of $\T^r$, and $\overline{\Gamma_{\theta}}$ simply is the cyclic subgroup $\langle \nu_{\theta} \rangle$, with Haar measure the uniform measure on it.\\

\begin{demo} We deal with the degenerate case and the non-degenerate case separately. For the degenerate case, we clearly have $\Gamma_{\theta} = \{(1, \dots, 1), \nu_{\theta}, \nu_{\theta}^2, \dots, \nu_{\theta}^{d-1}\} = \langle \nu_{\theta} \rangle$ so that if $f : \T^r \longrightarrow \C$ is any function (we don't even need it to be continuous), then we have \begin{align*}
\frac{1}{X} \sum_{n \leq X} f\left(e^{i\theta_1 n}, \dots, e^{i \theta_r n}\right) &= \frac{1}{X} \sum_{q=0}^{\left\lfloor \frac{X}{d}\right\rfloor} \sum_{a=0}^{d-1} f(\nu_{\theta}^{qd+a}) + o(1)\\ &=\frac{1}{X} \sum_{q=0}^{\left\lfloor \frac{X}{d}\right\rfloor} \sum_{a=0}^{d-1} f(\nu_{\theta}^a) + o(1)\\ &= \frac{\left\lfloor \frac{X}{d}\right\rfloor}{X} \sum_{a=0}^{d-1} f(\nu_{\theta}^a) + o(1) \underset{X \to +\infty}{\longrightarrow} \frac{1}{d} \sum_{a=0}^{d-1} f(\nu_{\theta}^a)
\end{align*} which is precisely what we needed to prove.\\

Now assume we are in the non-degenerate case, so that $m \geq 1$ and $\theta_1$ is not a rational multiple of $\pi$. We first show that $$\Gamma_{\theta} = \bigcup_{a=0}^{d-1} \nu_{\theta}^a \tilde{H}$$ where $$\tilde{H} := \left\{\left(e^{id\theta_1q}, \dots, e^{id\theta_mq}, \dots,  \prod_{k=1}^m e^{i h_{k, j} \theta_k q}, \dots \right) \mid q \in \Z \right\}.$$ To do this, we split $\Gamma_{\theta}$ according to its congruence classes modulo $d$ : $$\Gamma_{\theta} = \bigcup_{a=0}^{d-1} \Gamma_a,$$ where $$\Gamma_a := \left\{\left(e^{i\theta_1X}, \dots, e^{i\theta_rX}\right) \mid X \equiv a \text{ mod } d\right\}.$$ Expressing each $\theta_j$ with $m+1 \leq j \leq r$ in the basis $\{2\pi, \theta_1, \dots, \theta_m\}$ with \ref{coeff}, we find for $0 \leq a \leq d-1$, \begin{align*}\Gamma_a &= \left\{\left(e^{ia\theta_1}e^{id\theta_1q}, \dots, e^{ia\theta_m}e^{id\theta_mq}, \dots, e^{2i\pi l_j q} e^{2i\pi a c_j} \prod_{k=1}^m e^{iab_{k, j} \theta_k} \prod_{k=1}^m e^{i h_{k, j} \theta_k q}, \dots \right) \mid q \in \Z \right\}\\ &= \left\{\left(\nu_1^ae^{id\theta_1q}, \dots, \nu_m^ae^{id\theta_mq}, \dots, \nu_j^a \prod_{k=1}^m e^{i h_{k, j} \theta_k q}, \dots \right) \mid q \in \Z \right\} \end{align*} where we wrote $$\nu_j = e^{i \theta_j} \text{ for } 1 \leq j \leq m$$ and $$\nu_j = e^{2 i \pi c_j} \prod_{k=1}^m e^{ib_{k,j} \theta_k} = e^{i \theta_j} \text{ for } m+1 \leq j \leq r.$$ We have thus shown that $$\Gamma_{\theta} = \bigcup_{a=0}^{d-1} \nu_{\theta}^a \tilde{H}$$ as announced.

Let us show that this union is disjoint : let $a, b \in \{0, \dots, d-1\}$. If $\nu_1^a e^{id\theta_1 q} =  \nu_1^be^{id\theta_1 q'}$ then we have $\theta_1(a+qd-b-q'd) = 2k \pi$ for some $k \in \Z$. Since $\theta_1$ is not a rational multiple of $\pi$, we find $a+qd=b+q'd$, hence $a=b$ by uniqueness of the remainder in euclidean division.

Now the discrete Kronecker-Weyl theorem \ref{KW} implies that $\left\{\left(e^{i\theta_1X}, \dots, e^{i\theta_mX}\right) \mid X \in \Z\right\}$ is equidistributed in $\T^m$. Lifting by the continuous surjective homomorphism $$\begin{array}{rccl}
&\T^m &\longrightarrow &H_{\theta}\\\varphi : &(z_1, \dots, z_m)&\mapsto & \left(z_1^d, \dots, z_m^d, \prod_{k=1}^m z_k^{h_{k, m+1}}, \dots, \prod_{k=1}^m z_k^{h_{k, r}}\right)
\end{array}$$ we find that for every continuous $f : \T^r \rightarrow \C$, one has \begin{align*}\frac{1}{X} \sum_{q \leq X} f\left(e^{id\theta_1q}, \dots, e^{id\theta_mq}, \dots,  \prod_{k=1}^m e^{i h_{k, j} \theta_k q}, \dots \right) &= \frac{1}{X}\sum_{q \leq X} f \circ \varphi\left(e^{i\theta_1q}, \dots, e^{i\theta_mq}\right)\\ &\underset{X \to +\infty}{\longrightarrow} \int_{\T^m} f \circ \varphi \,\mathrm{d}\mu\\ &= \int_{H_{\theta}} f \,\mathrm{d}(\varphi_*\lambda).\end{align*} where $\varphi_*\lambda$ is the pushforward measure of the Lebesgue measure on $\T^m$ by $\varphi$. This measure is readily verified to be the Haar measure $\mu_{H_{\theta}}$ on $H_{\theta}$, since it has mass one and it is invariant by translations.

We have thus shown that for every continuous $f : \T^r \rightarrow \C$, one has $$\frac{1}{X} \sum_{q=1}^X f\left(e^{id\theta_1q}, \dots, e^{id\theta_mq}, \dots,  \prod_{k=1}^m e^{i h_{k, j} \theta_k q}, \dots \right) \underset{X \to +\infty}{\longrightarrow} \int_{H_{\theta}} f \,\mathrm{d}\mu_{H_{\theta}},$$ \textit{i.e.} that $\tilde{H}$ is equidistributed in $H_{\theta}$ with respect to its Haar measure. If we take any such $f$ and sum it over $\Gamma_{\theta}$ instead, we now find, using the previous disjoint decomposition of $\Gamma_{\theta}$, that \begin{align*}
&\frac{1}{X} \sum_{n \leq X} f\left(e^{i\theta_1n}, \dots, e^{i\theta_rn}\right)\\ &= \frac{1}{X} \sum_{a=0}^{d-1} \sum_{q=0}^{\left\lfloor\frac{X}{d}\right\rfloor}  f\left(\nu_1^ae^{id\theta_1q}, \dots, \nu_m^ae^{id\theta_mq}, \dots, \nu_j^a \prod_{k=1}^m e^{i h_{k, j} \theta_k q}, \dots \right) + o(1)\\ &= \frac{1}{d} \sum_{a=0}^{d-1} \frac{1}{\left\lfloor \frac{X}{d}\right\rfloor} \sum_{q=1}^{\left\lfloor \frac{X}{d}\right\rfloor} f_a\left(e^{id\theta_1q}, \dots, e^{id\theta_mq}, \dots, \prod_{k=1}^m e^{i h_{k, j} \theta_k q}, \dots \right) + o(1)
\end{align*} where $f_a : z \mapsto f(\nu_{\theta}^az)$ for $0 \leq a \leq d-1$. We finally obtain $$\frac{1}{X} \sum_{n \leq X} f\left(e^{i\theta_1n}, \dots, e^{i\theta_rn}\right) \underset{X \to +\infty}{\longrightarrow} \frac{1}{d} \sum_{a=0}^{d-1} \int_{H_{\theta}} f_a \,\mathrm{d}\mu_{H_{\theta}} = \int_{\overline{\Gamma_{\theta}}} f \,\mathrm{d}\theta.$$
\end{demo}

\begin{rmq} \begin{itemize}
\item[i)] In the non-degenerate case, the element $\nu_{\theta}$ of the above theorem is not of finite order in $\T^r$, since $\nu_1 = e^{i \theta_1}$ has infinite order in $\T$, but it has finite order dividing $d$ in $\T^r/H_{\theta}$. Therefore the product of the two subgroups $\langle \nu_{\theta} \rangle$ and $H_{\theta}$ is indeed $\bigcup_{a=0}^{d-1} \nu_{\theta}^aH_{\theta}$.
\item[ii)] If $\{\pi, \theta_1, \dots, \theta_r\}$ is $\Q$-linearly independent, then we have $m=r, d=1$ and $\overline{\Gamma_{\theta}} = \T^r$ as in Theorem \ref{KW}.
\item[iii)] The subgroup $H_{\theta}$ is a subtorus of $\T^r$ of dimension $m$, and the conclusion of the above theorem is that $\Gamma_{\theta}$ is equidistributed in the union of $d$ translates of this subtorus. Note that this union is not necessarily disjoint, as it would imply $\overline{\Gamma_{\theta}}$ has exactly $d$ connected components, but the number $d$ can be modified by choosing a different basis of $\Span_{\Q}(\pi, \theta_1, \dots, \theta_r)$ without changing $\Gamma_{\theta}$. In fact, it can easily be seen that in general, some of these translates may be equal. The exact number of connected components was already determined by Weyl (\cite[Satz 18]{Weyl}). Let $C=\{x \in \Q \mid \exists b \in \Z^r, \langle \theta, b \rangle = 2\pi x\}$ be the set of all $2\pi$ coefficients in the rational linear relations between $\theta_1, \dots, \theta_r$ and $2\pi$, after clearing denominators in the $\theta_i$ coefficients. Then the number of connected components $c$ of $\overline{\Gamma_{\theta}}$ is the lowest common multiple of the denominators of the elements of $C$. This could be established also with our method. As finding this number requires knowing every rational linear relations between $\theta_1, \dots, \theta_r$ and $2\pi$, or one such "minimal" relation, which seems unpractical when applying our results, we prefer working with the unoptimal number $d$ instead of the number $c$.
\end{itemize}
\end{rmq}

We now define a $\T^r$-valued random vector $Z_{\theta}$ associated with $\theta_1, \dots, \theta_r$. In the degenerate case, $Z_{\theta}$ simply is a uniform random variable on the cyclic subgroup $\langle \nu_{\theta} \rangle$. In the non-degenerate case, we define \begin{equation} \label{VA} Z_{\theta} := \left(Z_1^d, \dots, Z_m^d, \dots, \prod_{k=1}^m Z_k^{h_{k,j}}, \dots\right)\end{equation} where $Z_1, \dots, Z_m$ are independent uniform random variables on $\T$. Note that in each case, $Z_{\theta}$ has $\mu_{\H_{\theta}}$ for its probability distribution

\begin{cor} \label{Espe}For any continuous $f : \T^r \rightarrow \C$, we have $$\frac{1}{X} \sum_{n \leq X} f\left(e^{i\theta_1n}, \dots, e^{i\theta_rn}\right) \underset{X \to +\infty}{\longrightarrow} \frac{1}{d} \sum_{a=0}^{d-1} \E\left(f\left(\nu_{\theta}^aZ_{\theta}\right)\right).$$
\end{cor}

\begin{demo} This is just a reformulation of Theorem \ref{NKW}, where we observe that the distribution of the random vector $\nu_{\theta}^aZ_{\theta}$ is simply the measure $\mu_a$.
\end{demo}

In the context of prime number races over function fields, one approaches prime counting functions by functions of the form $t \mapsto c + \sum_{j=1}^r a_j e^{i \theta_j t} + \overline{a_j} e^{-i \theta_j t}$, with real $c$ and complex $a_j$. Note that those are in particular polynomials in $e^{i \theta_1 t}, e^{-i \theta_1 t}, \dots, e^{i \theta_r t}, e^{-i \theta_1 t}$, \textit{i.e.} Laurent polynomials in $e^{i \theta_1 t}, \dots, e^{i \theta_1 r}$, for which we prove the following key elementary lemma.

\begin{lem} \label{Atom}Let $f \in \C(X_1, \dots, X_r)$ with no pole in $\T^r$. Assume we are in the non-degenerate case. Then for $0 \leq a \leq d-1$, one has $\P(f(\nu_{\theta}^a Z_{\theta})=0) = 0$ if and only if there exists $n \equiv a \emph{ mod } d$ such that $f\left(e^{i\theta_1 n}, \dots, e^{i \theta_r n}\right) \neq 0$.
\end{lem}

\begin{demo} Recall that $\Gamma_a = \left\{\left(e^{i\theta_1X}, \dots, e^{i\theta_rX}\right) \mid X \equiv a \text{ mod } d\right\}$ is equidistributed in $\nu_{\theta}^a H_{\theta}$ from Theorem \ref{NKW} and that the distribution of $Z_{\theta}$ is precisely the Haar measure on $H_{\theta}$. Therefore, if $f\left(e^{i\theta_1 n}, \dots, e^{i \theta_r n}\right) = 0$ for all $n \equiv a \text{ mod } d$, then by continuity of $f$ and density of $\Gamma_a$ in $\nu_{\theta}^a H_{\theta}$, we have $\P(f(\nu_{\theta}^a Z_{\theta})=0)=1$.

We prove the converse statement by induction on $m$. If $m=1$, then the equation $$f\left(\nu_{\theta, 1}^a z^d, \nu_{\theta, 2}^a z^{h_1,2}, \dots, \nu_{\theta, r}^a z^{h_{1,r}}\right) = 0$$ in $z$ reduces, by clearing denominators, to a polynomial equation $P(z) = 0$ with one unknown. This equation is non-trivial because, writing $n=qd+a$ with $q \in \Z$, we have \begin{align*}
f\left(e^{i\theta_1 n}, \dots, e^{i \theta_r n}\right) &= f\left(e^{i a \theta_1} e^{iq\theta_1 d}, e^{i(qd+a)(2\pi c_2 + b_{1, 2} \theta_1)}, \dots, e^{i(qd+a)(2\pi c_r + b_{1, r} \theta_1)}\right)\\ &= f\left(\nu_{\theta, 1}^a e^{iq\theta_1 d}, \nu_{\theta_2}^a e^{iq \theta_1 h_{1, 2}}, \dots, \nu_{\theta_r}^a e^{iq \theta_1 h_{1, r}}\right) \neq 0
\end{align*} by hypothesis, so that $P\left(e^{iq\theta_1}\right) \neq 0$. Therefore, this equation has a finite number of solutions in $\C$, and in particular in $\T$. Since $Z_1$ is uniform on the circle, we certainly have $$\P(f(\nu_{\theta}^a Z_{\theta}) = 0) = \P\left(f\left(\nu_{\theta, 1}^a Z_1^d, \nu_{\theta, 2}^aZ_1^{h_1,2}, \dots, \nu_{\theta, r}^a Z_1^{h_{1,r}}\right) = 0\right)=0.$$ Now assume the result is true for $m-1 \in \N$. As before, by clearing denominators, the equation $$f\left(\nu_{\theta,1}^a z_1^d, \dots, z_m^d, \dots, \nu_{\theta,j}^a \prod_{k=1}^m z_k^{h_{k,j}}, \dots\right) = 0$$ is equivalent to a non-zero polynomial equation $P(z_1, \dots, z_m)=0$ with $m$ unknowns. Moreover, the set $F$ of all $z_m \in \T$ such that $P(X_1, \dots, X_{m-1}, z_m) = 0$ is finite, since it is the zero set of the $P(X_1, \dots, X_{m-1}, Y) \in \C[X_1, \dots, X_{m-1}][Y]$, which is non-zero because as above we have $P\left(e^{iq\theta_1}, \dots, e^{iq\theta_m}\right) \neq 0$ by hypothesis on $n=qd+a$. By the Fubini-Tonelli theorem, we find \begin{align*}
\P\left(f\left(\nu_{\theta}^a Z_{\theta}\right) = 0\right) &= \P\left(f\left(\nu_{\theta, 1}^a Z_1^d, \dots, \nu_{\theta, m}^a Z_m^d, \dots, \nu_{\theta, j}^a \prod_{k=1}^m Z_k^{h_{k,j}}, \dots\right) = 0\right)\\ &= \P(P(Z_1, \dots, Z_m)=0)\\ &= \int_{\T^m} \mathbf 1_{P^{-1}(\{0\})}(z) \,\mathrm{d}z\\ &= \int_{\T \setminus F} \left(\int_{\T^{m-1}} \mathbf{1}_{P( \cdot, z_m)^{-1}(\{0\})}(z_1, \dots, z_{m-1}) \,\mathrm{d}z_1 \dots \mathrm{d}z_{m-1}\right) \, \mathrm{d}z_m
\end{align*} The inner integral is zero by the induction hypothesis, so we conclude that $\P\left(f\left(\nu_{\theta}^a Z_{\theta}\right) = 0\right)=0.$
\end{demo}

We can now prove the following theorem which allows us to pass from continuous functions to indicator functions of subsets of $\R^D$ defined by strict inequalities between functions as in Lemma \ref{Atom}.

\begin{thm} \label{Prob} For $1 \leq j \leq D$, let $f_j \in \C(X_1, \dots, X_r)$ be real-valued and without pole on $\T^r$ and let $F_j : t \mapsto f_j\left(e^{i\theta_1 t}, \dots, e^{i\theta_rt}\right)$. Then we have $$\frac{1}{X} \#\left\{n \leq X \mid F_1(n) > \dots > F_D(n)\right\} \underset{X \to +\infty}{\longrightarrow} \frac{1}{d} \sum_{a=0}^{d-1} \P(f_1(\nu_{\theta}^aZ_{\theta}) > \dots > f_D(\nu_{\theta}^aZ_{\theta})).$$
\end{thm}

\begin{demo} In the degenerate case, we simply note that Corollary \ref{Espe} holds also for any function on $\T^r$, as in the beginning of the proof of Theorem \ref{NKW}, so we apply it to the function $\mathbf{1}_{x_1 > \dots > x_D} \circ (f_1, \dots, f_D)$.\\

Now assume we are in the non-degenerate case. We first remark that, by Lemma \ref{Atom} and its proof, if for some $0 \leq a \leq d-1$, there exists $1 \leq j \leq D-1$ such that $f_j\left(e^{in \theta_1}, \dots, e^{in \theta_r}\right) = f_{j+1}\left(e^{in \theta_1}, \dots, e^{in \theta_r}\right)$ for every $n \equiv a$ mod $d$, then $\P(f_j(\nu_{\theta}^a Z_{\theta}) = f_{j+1}(\nu_{\theta}^a Z_{\theta}))=1$, so that $$\P(f_1(\nu_{\theta}^aZ_{\theta}) > \dots > f_D(\nu_{\theta}^aZ_{\theta}))=0,$$ while we have $$\lim_{X \to +\infty} \frac{1}{X} \#\left\{n \leq X \mid n \equiv a\text{ mod } q, F_1(n) > \dots > F_D(n)\right\} = 0.$$
 Therefore, writing \begin{align*}&\frac{1}{X} \#\left\{n \leq X \mid F_1(n) > \dots > F_D(n)\right\}\\ &= \frac{1}{d} \sum_{a=0}^{d-1} \frac{1}{X} \#\left\{n \leq X \mid n \equiv a \text{ mod } d, F_1(n) > \dots > F_D(n)\right\} + o(1),\end{align*} and using the decomposition $\Gamma = \bigcup_{a=0}^{d-1} \nu_{\theta}^a \tilde{H}$ as in the proof of Theorem \ref{NKW}, we may assume that for every $0 \leq a \leq d-1, 1 \leq j \leq D-1$, there exists $n \equiv a$ mod $d$ such that $f_j\left(e^{in \theta_1}, \dots, e^{in \theta_r}\right) \neq f_{j+1}\left(e^{in \theta_1}, \dots, e^{in \theta_r}\right)$. By Lemma \ref{Atom} we then have $\P(f_j(\nu_{\theta}^a Z_{\theta}) = f_{j+1}(\nu_{\theta}^a Z_{\theta}))=0$ for every such $a$ and $j$.

Now, we need to approximate the indicator function $\mathbf 1_{x_1 > \dots > x_D}$ by continuous functions from above and below. We proceed in the following way : for every integer $k \geq 1$ and $x, y \in \R$, let $$g_k(x, y) := \left\{\begin{matrix} &0 \text{ if } x \leq y - \frac{1}{k}\\ &k(x-y)+1 \text{ if } y -\frac{1}{k} < x \leq y\\ &1 \text{ if } x > y.\end{matrix}\right.$$ Then for each integer $k \geq 1$, $g_k$ is continuous on $\R^2$ and for all $x, y \in \R$, $$\mathbf{1}_{x_1 > x_2}(x,y) \leq g_k(x, y) \leq \mathbf{1}_{x_1 > x_2 - \frac{1}{k}}(x, y).$$ For $k \geq 1$ let $G_k : (x_1, \dots, x_D) \mapsto \prod_{j=1}^{D-1} g_k(x_j, x_{j+1})$. Then for every $k \geq 1$ and $n \in \Z$, we have \begin{align*}\mathbf{1}_{x_1 > \dots > x_D}(F_1(n), \dots, F_D(n)) &= \prod_{j=1}^{D-1} \mathbf{1}_{x_j > x_{j+1}}(F_j(n), F_{j+1}(n))\\ &\leq \prod_{j=1}^{D-1} g_k(F_j(n), F_{j+1}(n))\\ &= G_k(F_1(n), \dots, F_D(n))\\ &\leq \prod_{j=1}^{D-1} \mathbf{1}_{x_j > x_{j+1} - \frac{1}{k}}(F_j(n), F_{j+1}(n))\\ &= \mathbf{1}_{x_1 > x_2 - \frac{1}{k} > \dots > x_D - \frac{D-1}{k}}(F_1(n), \dots, F_D(n)).\end{align*}

Now, by Corollary \ref{Espe}, for every $k \geq 1$, \begin{align*}&\limsup_{X \to +\infty} \frac{1}{X} \sum_{n=1}^X \mathbf{1}_{x_1 > \dots > x_D}\left(F_1(n), \dots, F_D(n)\right)\\ &\leq \limsup_{X \to +\infty} \frac{1}{X} \sum_{n=1}^X G_k(F_1(n), \dots, F_D(n))\\ &= \frac{1}{d} \sum_{a=0}^{d-1} \mathbb E(G_k(f_1(\nu_{\theta}^a Z_{\theta}), \dots, f_D(\nu_{\theta}^aZ_{\theta})))\\ &\leq \frac{1}{d} \sum_{a=0}^{d-1} \mathbb E\left( \mathbf{1}_{x_1 > x_2 - \frac{1}{k} > \dots > x_D - \frac{D-1}{k}}(f_1(\nu_{\theta}^a Z_{\theta}), \dots, f_D(\nu_{\theta}^aZ_{\theta})\right)\\&= \frac{1}{d} \sum_{a=0}^{d-1} \mathbb P\left(f_1(\nu_{\theta}^a Z_{\theta}) >  \dots > f_D(\nu_{\theta}^aZ_{\theta}) - \frac{D-1}{k}\right).\end{align*}

By downward continuity of $\P$, we get, by letting $k \to +\infty$, $$\limsup_{X \to +\infty} \frac{1}{X} \sum_{n=1}^X \mathbf{1}_{x_1 > \dots > x_D}\left(F_1(n), \dots, F_D(n)\right) \leq \frac{1}{d} \sum_{a=0}^{d-1} \mathbb P\left(f_1(\nu_{\theta}^a Z_{\theta}) \geq \dots \geq f_D(\nu_{\theta}^aZ_{\theta})\right).$$

Similarly, by considering the functions defined by $$(x, y) \mapsto \left\{\begin{matrix} &0 \text{ if } x < y\\ &k(x-y) \text{ if } y \leq x < y + \frac{1}{k}\\ &1 \text{ if } x > y + \frac{1}{k}\end{matrix}\right.,$$ we find $$\liminf_{X \to +\infty} \frac{1}{X} \sum_{n=1}^X \mathbf{1}_{x_1 > \dots > x_D}\left(F_1(n), \dots, F_D(n)\right) \geq \frac{1}{d} \sum_{a=0}^{d-1} \mathbb P\left(f_1(\nu_{\theta}^a Z_{\theta}) > \dots > f_D(\nu_{\theta}^aZ_{\theta})\right).$$

It remains to observe that the event $$\{f_1(\nu_{\theta}^a Z_{\theta}) \geq \dots \geq f_D(\nu_{\theta}^aZ_{\theta})\} \setminus \{f_1(\nu_{\theta}^a Z_{\theta}) > \dots > f_D(\nu_{\theta}^aZ_{\theta})\}$$ is included in $$\bigcup_{j=1}^{D-1} \{f_j(\nu_{\theta}^a Z_{\theta}) = f_{j+1}(\nu_{\theta}^a Z_{\theta})\}$$ which has probability zero, so that $$\mathbb P\left(f_1(\nu_{\theta}^a Z_{\theta}) \geq \dots \geq f_D(\nu_{\theta}^aZ_{\theta})\right) = \mathbb P\left(f_1(\nu_{\theta}^a Z_{\theta}) > \dots > f_D(\nu_{\theta}^aZ_{\theta})\right).$$
\end{demo}

In particular, we know that, with hypotheses from Theorem \ref{Prob}, the natural density $$\lim_{X \to +\infty} \frac{1}{X} \#\left\{n \leq X \mid F_1(n) > \dots > F_D(n)\right\}$$ always exists. We note however that this does not prove that every prime number race over a function field is weakly inclusive (as defined in Definition \ref{Incl}) because the normalized prime counting functions in this context have an extra $o(1)$ term (see Section \ref{2}). We will deal with such functions in Section \ref{Erreur}.

\subsection{The continuous case}

\label{Cont}

We now tackle the continuous case of the Kronecker-Weyl theorem.

Let $\theta_1, \dots, \theta_r$ be real numbers. Extract a basis $\{\theta_1, \dots, \theta_m\}$ of $\Span_{\Q}(\theta_1, \dots, \theta_r)$, and write \begin{equation} \label{coeffc} \theta_j = \sum_{k=1}^m b_{k,j} \theta_k \text{ for } m+1 \leq j \leq r\end{equation} with $b_{k,j} \in \Q$. Once again, we let \begin{equation} \label{dc} d := \mathrm{lcm}(\{\text{denominators of all } b_{k,j}\}),\end{equation} so that \begin{equation} \label{hc} h_{k,j} := db_{k,j} \in \Z.\end{equation}

The proof of the following theorem is similar to the proof of Theorem \ref{NKW}. We simply mention the necessary changes : discrete sums up to $X$ are replaced by integrals between $0$ and $X$, the splitting according to congruence classes modulo $d$ is replaced by the change of variable $y \longrightarrow dy$ and we appeal to the continuous version of the Kronecker-Weyl theorem (see \cite[Theorem 4.2]{Dev} or the remark at the end of the appendix). We note that most results in this section are made easier than in the previous section because linear dependence with $\pi$ doesn't have any effect on the continuous densities being studied.

\begin{thm} \label{CNKW} The one-parameter subgroup $$\Gamma_{\theta} = \left\{\left(e^{i\theta_1 y}, \dots, e^{i\theta_r y}\right) \mid y \in \R\right\}$$ is equidistributed in \begin{equation} \label{Hc} H_{\theta} = \left\{\left(z_1^d, \dots, z_m^d, \dots, \prod_{k=1}^m z_k^{h_{k,j}}, \dots\right) \mid (z_1, \dots, z_m) \in \T^m\right\},\end{equation} that is for every continuous $f : \T^r \longrightarrow \C$ one has $$\frac{1}{X} \int_0^X f\left(e^{i\theta_1 y}, \dots, e^{i \theta_r y}\right) \,\mathrm{d}y \underset{X \to +\infty}{\longrightarrow} \int_{H_{\theta}} f \,\mathrm{d}\mu_{H_{\theta}}$$ where $\mu_{H_{\theta}}$ is the normalized Haar measure on $H_{\theta}$.
\end{thm}

As in the discrete case, we define \begin{equation}\label{VAc}Z_{\theta} := \left(Z_1^d, \dots, Z_m^d, \dots, \prod_{k=1}^m Z_k^{h_{k,j}}, \dots\right)\end{equation} where $Z_1, \dots, Z_m$ are independent uniform random variables on $\T$, so that the Haar measure $\mu_{H_{\theta}}$ is the distribution of $Z_{\theta}$. We then obtain the following.

\begin{cor} \label{ConEsp} For any continuous $f : \T^r \rightarrow \C$, we have $$\frac{1}{X} \int_0^X f\left(e^{i\theta_1y}, \dots, e^{i\theta_ry}\right) \,\mathrm{d}y \underset{X \to +\infty}{\longrightarrow} \E\left(f\left(Z_{\theta}\right)\right).$$
\end{cor}


Note that the analog of Lemma \ref{Atom} holds with the only hypothesis that $f\left(e^{i \theta _1 y}, \dots, e^{i \theta_r y}\right) \neq 0$ for at least one $y \in \R$.

\begin{lem} \label{CAtom}Let $f \in \C(X_1, \dots, X_r)$ with no pole in $\T^r$. Then one has $\P(f(Z_{\theta})=0) = 0$ if and only if there exists $y \in \R$ such that $f\left(e^{i\theta_1 y}, \dots, e^{i \theta_r y}\right) \neq 0$.
\end{lem}

The proof of the continuous analog of Theorem \ref{Prob} then goes similarly. 

\begin{thm} \label{CProb} Let $f_1, \dots, f_D \in \C(X_1, \dots, X_r)$ be real-valued and without poles on $\T^r$ and let $F_j : t \mapsto f_j\left(e^{i\theta_1 t}, \dots, e^{i\theta_rt}\right)$. Then we have $$\frac{1}{X} \int_0^X \mathbf{1}_{x_1 > \dots > x_D}\left(F_1(y), \dots, F_D(y)\right) \,\mathrm{d}y \underset{X \to +\infty}{\longrightarrow} \P\left(f_1(Z_{\theta}) > \dots > f_D(Z_{\theta})\right).$$
\end{thm}

\subsection{The infinite-dimensional (continuous) case}

\label{Infini}

The goal of this section is to prove a result analogous to Theorem \ref{CProb} but with converging series in an infinite number of $e^{i \theta_n t}$, as a natural generalization of the above Laurent polynomials. This is the kind of functions we have to deal with in the context of prime number races over number fields, because the associated $L$-functions have an infinite number of non-trivial zeros. Those functions are often called \textit{almost periodic functions}.

\begin{defi} The $B^1$ semi-norm of a locally integrable function $f$ is $$||f||_{B^1} := \limsup_{X \to +\infty} \frac{1}{X} \int_0^X |f(y)| \,\mathrm{d}y.$$ A function $F : \R^+ \longrightarrow \R$ is said to be $B^1$-almost periodic if there exists a sequence $(P_N)_{N \geq 1}$ of trigonometric polynomials of the form $$P_N : t \mapsto \sum_{n=1}^{D_N} r_{n, N} e^{i \lambda_{n, N} t}$$ for some integer $D_N \geq 1$, $r_{n, N} \in \C$ and $\lambda_{n, N} \in \R$, such that $$||F-P_N||_{B^1} \underset{N \to +\infty}{\longrightarrow} 0.$$
\end{defi}

It turns out that prime counting functions over number fields are $B^2$-almost periodic, after applying the change of variable $x \to e^x$ (see \cite[Lemma 5.1.3]{Ng} or \cite[Proposition 4.4]{Dev} for a general statement), where $B^2$ semi-norm is defined by $||f||_{B^2} = \left(||f^2||_{B^1}\right)^{1/2}$. The Cauchy-Schwarz inequality easily implies that such functions are in particular $B^1$-almost periodic.\\

We simply quote the following important fact about $B^1$-almost periodic functions (\cite[p.104]{Bes}).

\begin{prop} Let $F : \R^+ \longrightarrow \R$ be $B^1$-almost periodic. There exists a countable set $\Lambda(F) = \{\lambda_n \mid n \geq 1\} \subset \R$ called the support of $F$, such that for every $n \geq 1$, $$a_n := \lim_{X \to +\infty} \frac{1}{X} \int_0^X f(y) e^{-i \lambda_n y} \,\mathrm{d}y \neq 0.$$ Moreover we have $||F - P_N||_{B^1} \underset{N \to +\infty}{\longrightarrow} 0,$ where $P_N(F) : t \mapsto \sum_{n=1}^N \left(a_n e^{-i \lambda_n t} + \overline{a_n} e^{-i \lambda_n t}\right).$
\end{prop}

The upshot of the above Proposition is that there exists a canonical way to approximate a given $B^1$-almost periodic function by trigonometric polynomials with respect to the $B^1$ semi-norm. We denote this fact by $F(t) \sim \sum_{n \geq 1} \left(a_n e^{-i \lambda_n t} + \overline{a_n} e^{-i \lambda_n t}\right)$. Note that this does not necessarily mean that the above series converges pointwise to $F$.\\

We begin by proving that $B^1$-almost periodic functions admit limiting distributions, using the language of random variables. The argument is essentially the one given in the proof of \cite[Theorem 2.9]{ANS}. We note that the left-hand side of formula $(2.10)$ in \textit{loc. cit.} should be replaced by $\limsup_{Y \to +\infty} \frac{1}{Y} \int_0^Y |\phi(y) - P_N(y)| \,\mathrm{d}y$, and that $Y$ should be assumed large enough in the last inequality in the proof of \cite[Theorem 2.9]{ANS}.

\begin{thm} \label{Helly} Let $F : \R^+ \longrightarrow \R$ be a $B^1$-almost periodic function. There exists a random variable $S$ such that for any continuous bounded function $g$ on $\R$ we have $$\frac{1}{X} \int_0^X g(F(y)) \,\mathrm{d}y \underset{X \to +\infty}{\longrightarrow} \E(g(S)).$$ In other words, $F$ admits $\mathbb P_S$, the distribution of $S$, as a limiting distribution.
\end{thm}

\begin{demo} The goal is to apply Corollary \ref{ConEsp} to each $P_N(F)$ before passing to the limit in $N$. Let $g$ be a bounded Lipschitz function on $\R$, with Lipschitz constant $c_g$. Then for any $N \geq 1$, $$\frac{1}{X} \int_0^X g(F(y)) \,\mathrm{d}y = \frac{1}{X} \int_0^X g(P_N(y)) \,\mathrm{d}y + \frac{1}{X} \int_0^X \left(g(F(y)) - g(P_N(y))\right) \,\mathrm{d}y.$$ By the triangular inequality one has $$\left|\frac{1}{X} \int_0^X \left(g(F(y)) - g(P_N(y))\right) \,\mathrm{d}y\right| \leq \frac{c_g}{X} \int_0^X |F(y) - P_N(y)| \,\mathrm{d}y$$ so that $\limsup_{X \to +\infty} \left|\frac{1}{X} \int_0^X \left(g(F(y)) - g(P_N(y))\right) \,\mathrm{d}y\right| \underset{N \to +\infty}{\longrightarrow} 0$. On the other hand, Corollary \ref{ConEsp} yields $$\frac{1}{X} \int_0^X g(P_N(y)) \,\mathrm{d}y \underset{X \to +\infty}{\longrightarrow} \E(g(S_N))$$ for some random variable $S_N$ built from the linear relations over $\Q$ between the real numbers $\lambda_{1, N}, \dots, \lambda_{D_N, N}$.

This proves that $$\limsup_{X \to +\infty} \frac{1}{X} \int_0^X g(F(y)) \,\mathrm{d}y = \E(g(S_N)) + o(1),$$ and similarly we have $$\liminf_{X \to +\infty} \frac{1}{X} \int_0^X g(F(y)) \,\mathrm{d}y = \E(g(S_N)) + o(1),$$ where both $o(1)$ are quantities going to zero as $N$ tends to infinity. Therefore $$\limsup_{X \to +\infty} \frac{1}{X} \int_0^X g(F(y)) - \liminf_{X \to +\infty} \frac{1}{X} \int_0^X g(F(y)) = 0$$ since it is independent of $N$ and $o(1)$ with respect to $N$. We have thus shown that the quantity $\frac{1}{X} \int_0^X g(F(y)) \,\mathrm{d}y$ admits a limit as $X$ tends to infinity, and $\E(g(S_N))$ converges to this (same) limit as $N$ tends to infinity.

We now prove that the sequence $(S_N)_{N \geq 1}$ converges in distribution to some random variable $Z$. To do so, we apply Prohorov's theorem \cite[Theorem 5.1]{Bill} (or Helly's selection theorem as it is called in \cite[Lemma 2.8]{ANS}), which in particular states that a tight sequence of probability measures on $\R$ admits a weakly converging subsequence. Recall that a family $(\mu_n)_{n \geq 1}$ of probability measures on $\R$ is tight when there is no "escape of mass to infinity" along the family, \textit{i.e.} for every $\varepsilon > 0$, one can find a compact $K \subset \R$ such that $\mu_n(K) \geq 1 - \varepsilon$ for every $n \geq 1$. Assuming that $\left(\mathbb P_{S_N}\right)_{N \geq 1}$ is tight, and denoting by $\mu$ the weak limit of one of its subsequence, $\left(\E(g(S_N))\right)_{N \geq 1}$ can only converge to $\int_{\R} g \,\mathrm{d}\mu$ when $g$ is a bounded continuous function on $\R$. Since this holds for every bounded Lipschitz function on $\R$ by the above computations, the Portmanteau theorem \cite[Theorem 2.1]{Bill} implies that $\left(\mathbb{P}_{S_N}\right)_{n \geq 1}$ converges weakly to $\mu$. Finally, the limit probability measure $\mu$ is the distribution of $S := F^{-1}(U)$, where $F$ is the distribution function of $\mu$, $F^{-1}$ its generalized inverse and $U$ is uniform on $[0, 1]$ (see \cite[Theorem 2.1]{Devr}), so that $\left(S_N\right)_{N \geq 1}$ converges in distribution to $S$.

It only remains to prove that $\left(\mathbb P_{S_N}\right)_{N \geq 1}$ is tight. Let $A > 0$. As a straightforward application of Theorem \ref{CProb}, we obtain $$\P(|S_N| > A) = \lim_{X \to +\infty} \frac{1}{X} \int_0^X \mathbf{1}_{x > A}(|P_N(y)|) \,\mathrm{d}y.$$ By Markov's inequality, we have for every $N \geq 1$ and $X > 0$, $$\frac{1}{X} \int_0^X \mathbf{1}_{x > A}(|P_N(y)|) \,\mathrm{d}y \leq \frac{1}{AX} \int_0^X |P_N(y)| \,\mathrm{d}y.$$ For every $N \geq 1$ and $y \in \R^+$, one has $|P_N(y)| \leq |F(y)| + |F(y) - P_N(y)|$. Now $L := \limsup_{X \to +\infty} \frac{1}{X} \int_0^X |F(y)| \,\mathrm{d}y < +\infty$ since $\limsup_{X \to +\infty} \frac{1}{X} \int_0^X |F(y) - P_n(y)| \,\mathrm{d}y < +\infty$ for at least one $n$, and $\limsup_{X \to +\infty} \frac{1}{X} \int_0^X |P_n(y)| \,\mathrm{d}y < +\infty$ by Corollary \ref{ConEsp}. Finally we obtain $$\P(|S_N| > A) \leq \frac{1}{A} \limsup_{X \to +\infty} \frac{1}{X} \int_0^X |F(y) - P_N(y)| \,\mathrm{d}y + \frac{L}{A} \ll \frac{1}{A}$$ independently of $N$, which proves the tightness of $\left(\mathbb P_{S_N}\right)_{N \geq 1}$.
\end{demo}

\begin{rmq} 
\begin{itemize}
\item[i)] The proof goes similarly for vector-valued $B^1$-almost periodic functions, as in \cite[Theorem 2.9]{ANS}.

\item[ii)]The key argument in the above proof was Prohorov's theorem, or Helly's selection theorem, but this is an indirect argument. Using our explicit version of the Kronecker-Weyl theorem, and with additional hypotheses on the sequence $(\theta_n)_{n \geq 1}$, we can give a more explicit description of the random variable $S$ in terms of the function $F$ (see Section \ref{Momente}).
\end{itemize}
\end{rmq}

The next step in our analysis is to pass from bounded continuous functions to indicator functions of sets defined by strict inequalities. Mimicking the proof of Theorem \ref{Prob}, we obtain the following.

\begin{prop} \label{Infineg} Let $F_1, \dots, F_D : \R^+ \longrightarrow \R$ be $B^1$-almost periodic functions with $\Lambda(F_j) \subset \Theta$ for $1 \leq j \leq D$. Let $S_1, \dots, S_D$ be the random variables associated to $F_1, \dots, F_D$ in Theorem \ref{Helly}. Then \begin{align*}\mathbb{P}(S_1 > \dots > S_D) &\leq \liminf_{X \to +\infty} \frac{1}{X} \int_0^X \mathbf{1}_{x_1 > \dots > x_D}(F_1(y), \dots, F_D(y)) \,\mathrm{d}y\\ &\leq \limsup_{X \to +\infty} \frac{1}{X} \int_0^X \mathbf{1}_{x_1 > \dots > x_D}(F_1(y), \dots, F_D(y)) \,\mathrm{d}y \leq \mathbb{P}(S_1 \geq \dots \geq S_D).\end{align*}
\end{prop}

\begin{rmq}
We cannot expect an equality without any hypothesis on $\Theta$. For instance, consider the case $D=2$, $F_2=0$ and $F_1$ a non-zero continuous function with compact support on $\R$. For a less trivial example, we can use an everywhere converging Fourier series which has a non-constant sum, but such that its sum is constant on a non-empty interval.\\
\end{rmq}

We now look for conditions on $\Theta$ to imply equality in the previous Proposition. Such a condition was found by Devin in \cite{Dev2} : if $\Span_{\Q} \Theta$ decomposes as a direct sum $\Span_{\Q} T \oplus \Span_{\Q}(\Theta \setminus T)$, where $T$ is a finite subset of $\Theta$, then the random variable associated to a non-constant $B^1$-almost periodic function $S$ as in Theorem \ref{Helly} does not admit atoms. Note that this hypothesis on $\Span_{\Q} \Theta$ is a weakening of the hypothesis of the existence of "self-sufficient zeros" in \cite{MN}. The proof of Devin consists in showing that the characteristic function of $S$ is decaying sufficiently fast at infinity, by using known bounds on oscillatory integrals, and then using a lemma of Wiener, relating this decay to the continuity of the distribution of $S$ (see the proof of Theorem 1.2 and Corollary 1.4 in \cite{Dev2}). Our method allows us to show the same thing but with a considerably simpler proof thanks to Lemma \ref{CAtom}.

\begin{thm}[Devin] \label{Conv} Assume $\Span_{\Q}(\Theta) = \Span_{\Q} (\theta_1, \dots, \theta_m) \oplus \Span_{\Q}(\{\theta_n \mid n > m\})$. Let $F : \R^+ \longrightarrow \R$ be such that $||F - P_N||_{B^1} \underset{N \to +\infty}{\longrightarrow} 0$ where $$P_N : t \mapsto c + \sum_{n \leq N} \left(a_n e^{i \theta_n t} + \overline{a_n} e^{-i \theta_n t}\right).$$ Let $\mu_F, \mu_{F-P_m}$ and $\mu_{P_m}$ be the limiting distributions of $F, F-P_m$ and $P_m$ respectively. Then $\mu_F = \mu_{F-P_m} * \mu_{P_m}$. In particular, if $P_m$ is not constant, then $\mu_F(\{x\}) = 0$ for any $x \in \R$.
\end{thm}

\begin{demo} The function $F$ and $R=F - P_m$ are $B^1$-almost periodic functions, and therefore admit limiting distributions $\mu_F$ and $\mu_{R}$ by Theorem \ref{Helly}. On the other hand, $P_m$ also admits a limiting distribution because of Theorem \ref{CNKW}, say $\mu_{P_m}$.

For $N > m$, let $P_N' : t \mapsto \sum_{m < n \leq N} a_n e^{i \theta_n t} + \overline{a_n} e^{-i \theta_n t}$. Just as $P_m$, those admit limiting distributions $\mu_{P_N'}$ and by the proof of Theorem \ref{Helly}, $\left(\mu_{P_N'}\right)_{N > m}$ converges weakly to $\mu_R$. Also, for every $N > m$, $\mu_{P_N' + P_m} = \mu_{P_N'} * \mu_{P_m}$ because, by construction and the hypothesis on $\Span_{\Q} \Theta$, every $Z_{n, N}, 1 \leq n \leq m$ is independent of every $Z_{n', N}, m < n' \leq N$. As above, $\left(\mu_{P_N' + P_m}\right)_{N > m}$ converges weakly to $\mu_F$.

For any Borel probability measure on $\R$, let $\hat{\mu}$ be its characteristic function. Then for every $N > m, \widehat{\mu_{P_N'} * \mu_{P_m}} = \hat{\mu}_{P_N'} \hat{\mu}_{P_m}$ and this converges pointwise to $\hat{\mu}_{R} \hat{\mu}_{P_m} = \widehat{\mu_{R} * \mu_{P_m}}$. By Lévy's continuity theorem (\cite[Theorem 3.3.17]{Dur}), this means that $(\mu_{P_N'} * \mu_{P_m})_{N > m}$ converges weakly to $\mu_{R} * \mu_{P_m}$. Therefore we have proved that $\mu_F = \mu_R * \mu_{P_m}$.

Now if $P_m$ is not constant, then by Lemma \ref{CAtom}, we have $\mu_{P_m}(\{y\})=0$ for every $y \in \R$, and thus for any $x \in \R$, $$\mu_F(\{x\}) = (\mu_{R} * \mu_{P_m})(\{x\}) = \int_{\R} \mu_{P_m}(\{x-y\}) \,\mathrm{d}\mu_{R}(y) = 0.$$
\end{demo}


Combining the previous Theorem with Proposition \ref{Infineg} we obtain the following.

\begin{cor} Assume $\Span_{\Q} \Theta = \Span_{\Q} T \oplus \Span_{\Q}(\Theta \setminus T)$ for some non-empty finite subset $T$ of $\Theta$. Let $F_1, \dots, F_D : \R^+ \longrightarrow \R$ be $B^1$-almost periodic functions such that for $1 \leq j \leq D, F_j \sim c_j + \sum_{\theta \in \Theta} \left(a_{\theta, j} e^{i \theta t} + \overline{a_{\theta, j}} e^{-i \theta t}\right).$ Let $S_1, \dots, S_D$ be the random variables associated to $F_1, \dots, F_D$ in Theorem \ref{Helly}. If for every $1 \leq j \leq D-1$,  the function $t \mapsto \sum_{\theta \in T} (a_{\theta, j}-a_{\theta, j+1}) e^{i \theta y} + (\overline{a_{\theta, j}} - \overline{a_{\theta, j+1}}) e^{-i \theta y}$ is not constant, then $$\lim_{X \to +\infty} \frac{1}{X} \int_0^X \mathbf{1}_{x_1 > \dots > x_D}(F_1(y), \dots, F_D(y)) \,\mathrm{d}y = \mathbb{P}(S_1 > \dots > S_D).$$
\end{cor}

The condition that a certain linear combination of $e^{i \theta t}$ is not constant can be easily translated by the non-vanishing of its coefficients.

\begin{lem} Let $T$ be a finite subset of $\Theta$ and for each $\theta \in T$, let $a_{\theta} \in \C$. If the function $t \mapsto \sum_{\theta \in T} \left(a_{\theta} e^{i \theta t} + \overline{a_{\theta}} e^{-i \theta t}\right)$ is constant then $a_{\theta} = 0$ for every $\theta \in T$.
\end{lem}

\begin{demo} Let $P : t \mapsto \sum_{\theta \in T} a_{\theta} e^{i \theta t} + \overline{a_{\theta}} e^{-i \theta t}$ and assume it is constant. Then all of its derivatives are zero. In particular, we have for $1 \leq k \leq 2\#T$, $$P^{(k)}(0) = \sum_{\theta \in T} (a_{\theta} + (-1)^k \overline{a_{\theta}}) \theta^k = 0,$$ which means that $(a_{\theta} + \overline{a_{\theta}})_{\theta \in T}$ and $(a_{\theta} - \overline{a_{\theta}})_{\theta \in T}$ are both solutions of Vandermonde linear systems with non-zero determinant since the elements of $\Theta$, and therefore of $T$, are non-zero and pairwise distinct. This implies $a_{\theta} = \overline{a_{\theta}}$ and $a_{\theta} = - \overline{a_{\theta}}$ for every $\theta \in T$, and thus $a_{\theta} = 0$ for every $\theta \in T$.
\end{demo}

\begin{rmq} This lemma can also be seen as an application of Artin's lemma on the linear independence of characters (see \cite[VI, Theorem 4.1]{Lang}).
\end{rmq}

\begin{cor} \label{Probinf} Assume $\Span_{\Q} \Theta = \Span_{\Q} T \oplus \Span_{\Q}(\Theta \setminus T)$ for some non-empty finite subset $T$ of $\Theta$. Let $F_1, \dots, F_D : \R^+ \longrightarrow \R$ be $B^1$-almost periodic functions such that for $1 \leq j \leq D, F_j \sim c_j + \sum_{\theta \in \Theta} \left(a_{\theta, j} e^{i \theta t} + \overline{a_{\theta, j}} e^{-i \theta t}\right).$ Let $S_1, \dots, S_D$ be the random variables associated to $F_1, \dots, F_D$ in Theorem \ref{Helly}. If for every $1 \leq j \leq D-1$, there exists $\theta \in T$ such that $a_{\theta, j} \neq a_{\theta, j+1}$ then $$\lim_{X \to +\infty} \frac{1}{X} \int_0^X \mathbf{1}_{x_1 > \dots > x_D}(F_1(y), \dots, F_D(y)) \,\mathrm{d}y = \mathbb{P}(S_1 > \dots > S_D).$$
\end{cor}

\begin{demo} Simply combine the previous two results.
\end{demo}

\section{Applications}

\subsection{Quantities with an error term}

\label{Erreur}

We now investigate to what extent Theorem \ref{Prob}, Theorem \ref{CProb} and Corollary \ref{Probinf} still hold for functions with an extra error term, because the explicit formulas for prime counting functions admit such an error term. To simplify notations, when $G_1, \dots, G_D$ are functions and $\mathcal{R}$ is a $D$-ary relation ($\mathcal R$ will be $\{(x_1, \dots, x_D) \in \R^D \mid x_1 > \dots > x_D\}$ or $\{(x_1, \dots, x_D) \in \R^D \mid x_1 \geq \dots \geq x_D\}$ below), we set $$\underline{\delta}(\mathcal{R}(G_1, \dots, G_D)) := \liminf_{X \to +\infty} \frac{1}{X} \#\{n \in \{1, \dots, X\} \mid \mathcal{R}(G_1(n), \dots, G_D(n))\}$$ and $$\overline{\delta}(\mathcal{R}(G_1, \dots, G_D)) := \limsup_{X \to +\infty} \frac{1}{X} \#\{n \in \{1, \dots, X\} \mid \mathcal{R}(G_1(n), \dots, G_D(n))\}$$ When $\underline{\delta}(\mathcal{R}(G_1, \dots, G_D)) = \overline{\delta}(\mathcal{R}(G_1, \dots, G_D))$ we denote by $\delta(\mathcal{R}(G_1, \dots, G_D))$ their common value.

We start with the discrete case.

\begin{thm} \label{Reste} Let $\theta_1, \dots, \theta_r$ be real numbers and $f_1, \dots, f_D \in \C(X_1, \dots, X_r)$ be real-valued and without pole on $\T^r$. Let $G : t \mapsto (F_1(t), \dots, F_D(t)) + r(t)$ where $F_j(t) = f_j\left(e^{i \theta_1 t}, \dots, e^{i \theta_r t}\right)$ for $1 \leq j \leq D$, and $r(t) = (r_1(t), \dots, r_D(t))=o(1)$ as $t \to +\infty$. 
\begin{itemize}
\item[i)] Degenerate case : Assume that $\theta_i \in \pi\Q$ for $1 \leq i \leq r$. Then 
\begin{align*}\mathbb P(f_1(Z_{\theta}) > \dots > f_D(Z_{\theta})) &\leq \underline{\delta}(G_1 > \dots > G_D)\\ &\leq \overline{\delta}(G_1 > \dots > G_D) \leq \mathbb P(f_1(Z_{\theta}) \geq \dots \geq f_D(Z_{\theta})).\end{align*} In particular, if for every $1 \leq j \leq D-1$ and every $n \in \Z$ one has $F_j(n) \neq F_{j+1}(n)$, then $\delta(G_1 > \dots > G_D)$ exists and we have $$\delta(G_1 > \dots > G_D) = \mathbb P(f_1(Z_{\theta}) > \dots > f_D(Z_{\theta})).$$

\item[ii)] Non-degenerate case : Assume $\theta_i \not \in \pi \Q$ for at least one $i \in \{1, \dots, r\}$. Then 
\begin{align*}\frac{1}{d} \sum_{a=0}^{d-1} \mathbb P(f_1(\nu_{\theta}^aZ_{\theta}) > \dots > f_D(\nu_{\theta}^a Z_{\theta})) &\leq \underline{\delta}(G_1 > \dots > G_D)\\ &\leq \overline{\delta}(G_1 > \dots > G_D)\\ &\leq \frac{1}{d} \sum_{a=0}^{d-1} \mathbb P(f_1(\nu_{\theta}^a Z_{\theta}) \geq \dots \geq f_D(\nu_{\theta}^a Z_{\theta})).\end{align*}

Moreover, if for every $1 \leq j \leq D-1$ and $0 \leq a \leq d-1$, there exists $n \equiv a \emph{ mod } d$ such that $F_j(n) \neq F_{j+1}(n)$, then $\delta(G_1 > \dots > G_D)$ exists and $$\delta(G_1 > \dots > G_D) = \frac{1}{d} \sum_{a=0}^{d-1} \mathbb P(f_1(\nu_{\theta}^a Z_{\theta}) > \dots > f_D(\nu_{\theta}^a Z_{\theta})).$$
\end{itemize}
\end{thm}

\begin{demo} Let $\varepsilon > 0$. There exists $n_0 \geq 1$ such that for every $n \geq n_0$ and $1 \leq j \leq D$, we have $|r_j(n)| < \varepsilon$. Now for every $n \geq n_0$, the inequalities $$F_1(n) > F_2(n) + 2\varepsilon > \dots > F_D(n) + 2(D-1) \varepsilon$$ imply $$G_1(n) > G_2(n) > \dots > G_D(n)$$ which in turn imply $$F_1(n) > F_2(n) - 2\varepsilon > \dots > F_D(n) - 2(D-1) \varepsilon.$$

In the degenerate case, we get, using Theorem \ref{Prob}, that \begin{align*}&\mathbb P(f_1(Z_{\theta}) > \dots > f_D(Z_{\theta}) + 2(D-1)\varepsilon)\\ &\leq \underline{\delta}(G_1 > \dots > G_D)\\ &\leq \overline{\delta}(G_1 > \dots > G_D) \leq \mathbb P(f_1(Z_{\theta}) > \dots > f_D(Z_{\theta})-2(D-1)\varepsilon),\end{align*} while in the non-degenerate case, we get, still using Theorem \ref{Prob}, that \begin{align*}&\frac{1}{d} \sum_{a=0}^{d-1} \mathbb P(f_1(\nu_{\theta}^aZ_{\theta}) > \dots > f_D(\nu_{\theta}^a Z_{\theta}) + 2(D-1)\varepsilon)\\ &\leq \underline{\delta}(G_1 > \dots > G_D)\\ &\leq \overline{\delta}(G_1 > \dots > G_D)\\ &\leq \frac{1}{d} \sum_{a=0}^{d-1} \mathbb P(f_1(\nu_{\theta}^a Z_{\theta}) > \dots > f_D(\nu_{\theta}^a Z_{\theta}) - 2(D-1)\varepsilon).\end{align*} In both cases, we obtain the announced inequalities on $\underline{\delta}(G_1 > \dots > G_D)$ and $\overline{\delta}(G_1 > \dots > G_D)$ by letting $\varepsilon$ tend to $0$ as in the proof of Theorem \ref{Prob}.

Finally, the last hypotheses imply that $\P(f_1(Z_{\theta}) > \dots > f_D(Z_{\theta})) = \P(f_1(Z_{\theta}) \geq \dots \geq f_D(Z_{\theta}))$ in the degenerate case since $Z_{\theta}$ is uniform on $\langle \nu_{\theta} \rangle$, while they imply $\P(f_1(\nu_{\theta}^a Z_{\theta}) > \dots > f_D(\nu_{\theta}^aZ_{\theta})) = \P(f_1(\nu_{\theta}^a Z_{\theta}) \geq \dots \geq f_D(\nu_{\theta}^a Z_{\theta}))$ for every $0 \leq a \leq D-1$ in the non-degenerate case because of Lemma \ref{Atom}.
\end{demo}

The proof of the next theorem is completely similar, based on Theorem \ref{CProb} and Corollary \ref{Probinf}. This time, $\delta(\mathcal{R}(G_1, \dots, G_D))$ means, when it exists, $$\lim_{X \to +\infty} \frac{1}{X} \int_0^X \mathbf{1}_{\mathcal{R}}(G_1(y), \dots, G_D(y)) \,\mathrm{d}y,$$ and $\underline{\delta}(\mathcal{R}(G_1, \dots, G_D))$ and $\overline{\delta}(\mathcal{R}(G_1, \dots, G_D))$ the corresponding $\liminf$ and $\limsup$.

\begin{thm} \label{ResteC}
\begin{itemize}
\item[i)] Let $\theta_1, \dots, \theta_m$ be real numbers, $f_1, \dots, f_D \in \C(X_1, \dots, X_r, )$ be pairwise distinct and real-valued on $\T^r$. Let $G : t \mapsto (F_1(t), \dots, F_D(t)) + o(1)$ as $t \to +\infty$, where $F_j(t) = f_j\left(e^{i \theta_1 t}, \dots, e^{i \theta_r t}\right)$ for $1 \leq j \leq D$. Then 
\begin{align*}\mathbb P(f_1(Z_{\theta}) > \dots > f_D(Z_{\theta})) &\leq \underline{\delta}(G_1 > \dots > G_D)\\ &\leq \overline{\delta}(G_1 > \dots > G_D) \leq \mathbb P(f_1(Z_{\theta}) \geq \dots \geq f_D(Z_{\theta})).\end{align*} Moreover, if for every $1 \leq j \leq D-1$, there exists $y \in \R$ such that $F_j(y) \neq F_{j+1}(y)$, then $\delta(G_1 > \dots > G_D)$ exists and we have $$\delta(G_1 > \dots > G_D) = \mathbb P(f_1(Z_{\theta}) > \dots > f_D(Z_{\theta})).$$

\item[ii)] Let $\theta = (\theta_n)_{n \geq 1}$ be a sequence of pairwise distinct positive real numbers and $\Theta = \{\theta_n \mid n \geq 1\}$. Let $F_1, \dots, F_D : \R^+ \longrightarrow \R$ be $B^1$-almost periodic functions such that for $1 \leq j \leq D, F_j \sim c_j + \sum_{\theta \in \Theta} \left(a_{\theta, j} e^{i \theta t} + \overline{a_{\theta, j}} e^{-i \theta t}\right)$ and $G : t \mapsto (F_1(t), \dots, F_D(t)) + o(1)$ as $t \to +\infty$. Then, with $S_1, \dots, S_D$ the random variables associated to $F_1, \dots, F_D$ in Theorem \ref{Helly}, we have \begin{align*}\mathbb{P}(S_1 > \dots > S_D) &\leq \underline{\delta}(G_1 > \dots > G_D)\\ &\leq \overline{\delta}(G_1 > \dots > G_D) \leq \mathbb{P}(S_1 \geq \dots \geq S_D).\end{align*}

Moreover, if $\Span_{\Q} \Theta = \Span_{\Q} T \oplus \Span_{\Q}(\Theta \setminus T)$ for some non-empty finite subset $T$ of $\Theta$, and if for every $1 \leq j \leq D-1$, there exists $\theta \in T$ such that $a_{\theta, j} \neq a_{\theta, j+1}$, then $\delta(G_1 > \dots > G_D)$ exists and $$\delta(G_1 > \dots > G_D) = \P(S_1 > \dots > S_D).$$
\end{itemize}
\end{thm}

\subsection{Non-critical densities}

\label{Noncrit}

We give sufficient conditions for the densities we study to be positive.

\begin{prop} \label{PosDisc} Let $\theta_1, \dots, \theta_r$ be real numbers. Let $f_1, \dots, f_D, F_1, \dots, F_D$ and $G_1, \dots, G_D$ be as in Theorem \ref{Reste}.

\begin{itemize}
\item[i)] Degenerate case : Assume that $\theta_i \in \pi\Q$ for $1 \leq i \leq r$. If there exists $n \in \Z$ such that $F_1(n) > \dots > F_D(n)$, then $0 < \underline{\delta}(G_1 > \dots > G_D)$ and if there exists $n \in \Z$ such that $F_1(n) \geq \dots \geq F_D(n)$ does not hold, then $\overline{\delta}(G_1 > \dots > G_D) < 1$.

\item[ii)] Non-degenerate case : Assume $\theta_i \not \in \pi \Q$ for at least one $i \in \{1, \dots, r\}$. If there exist $a \in \{0, \dots, d-1\}$ and $z \in H_{\theta}$ such that $f_1(\nu_{\theta}^az) > \dots > f_D(\nu_{\theta}^az)$ then $0 < \underline{\delta}(G_1 > \dots > G_D)$. Also, if there exist $a \in \{0, \dots, d-1\}$ and $z \in H_{\theta}$ such that $f_1(\nu_{\theta}^az) \geq \dots \geq f_D(\nu_{\theta}^az)$ does not hold, then $\overline{\delta}(G_1 > \dots > G_D) > 1$. In particular, if there exists $n \in \Z$ such that $F_1(n) > \dots > F_D(n)$, then $0 < \underline{\delta}(G_1 > \dots > G_D)$, and if there exists $n \in \Z$ such that $F_1(n) \geq \dots \geq F_D(n)$ does not hold, then $\overline{\delta}(G_1 > \dots > G_D) < 1$.
\end{itemize}

\end{prop}

\begin{demo}\begin{itemize}
\item[i)] It is immediate by Theorem \ref{Prob} and Theorem \ref{Reste} i).

\item[ii)] Assume that there exist $a \in \{0, \dots, d-1\}$ and $z \in H_{\theta}$ such that $f_1(\nu_{\theta}^az) > \dots > f_D(\nu_{\theta}^az)$. By continuity of $f_1, \dots, f_D$, there exists an open subset $U$ of $H_{\theta}$ such that $z \in U$ and for all $z' \in U, f_1(\nu_{\theta}^az') > \dots > f_D(\nu_{\theta}^az')$. Therefore, we have \begin{align*}
\mathbb{P}(f_1(\nu_{\theta}^aZ_{\theta}) > \dots > f_D(\nu_{\theta}^aZ_{\theta})) &\geq \mathbb{P}(Z_{\theta} \in U)\\ &= \mathbb{P}((Z_1, \dots, Z_m) \in \varphi^{-1}(U))\\ &=\lambda(\varphi^{-1}(U)) > 0,
\end{align*} where $\varphi$ is the continuous map $$\begin{array}{rccl}
&\T^m &\longrightarrow &H_{\theta}\\\varphi : &(z_1, \dots, z_m)&\mapsto & \left(z_1^d, \dots, z_m^d, \prod_{k=1}^m z_k^{h_{k, m+1}}, \dots, \prod_{k=1}^m z_k^{h_{k, r}}\right)
\end{array}$$ and $\lambda$ is the Lebesgue measure on $\T^m$. By Theorem \ref{Reste} \textit{ii)}, we obtain $\underline{\delta}(F_1 > \dots > F_D) \geq \frac{\mathbb{P}(f_1(\nu_{\theta}^aZ_{\theta}) > \dots > f_D(\nu_{\theta}^aZ_{\theta}))}{d} > 0.$ The proof of the second statement is similar since the negation of $f_1(\nu_{\theta}^az) \geq \dots \geq f_D(\nu_{\theta}^az)$ is also an open condition on $z$.

The last stament is immediate since, if $n \equiv a \text{ mod } d$, then $\left(e^{i\theta_1 n}, \dots, e^{i \theta_rn}\right) = \nu_{\theta}^a z$ for some $z \in H_{\theta}$ by construction.
\end{itemize}
\end{demo}

\begin{rmq} \begin{itemize}
\item[i)] In the degenerate case, we actually have the lower bound $\frac{1}{d} \leq \underline{\delta}(G_1 > \dots > G_D)$ whenever $\underline{\delta}(G_1 > \dots > G_D) > 0$, and the upper bound $\overline{\delta}(G_1 > \dots > G_D) \leq 1 - \frac{1}{d}$ whenever $\overline{\delta}(G_1 > \dots > G_D) < 1$.

\item[ii)] The converses of the above statements are false in general. For example it may happen that $f_1(z) = \dots = f_D(z)$ for every $z \in \T^r$, but $\underline{\delta}(G_1 > \dots > G_D) > 0$ because $r_1(n) > \dots > r_D(n)$ for a positive proportion of $n \in \N$.

\end{itemize}
\end{rmq}

The continuous version of the above Proposition is proved in the same way as \textit{ii)} above.

\begin{prop} \label{PosCont} Let $\theta_1, \dots, \theta_r$ be real numbers. Let $f_1, \dots, f_D, F_1, \dots, F_D$ and $G$ be as in Theorem \ref{ResteC} i). If there exists $z \in H_{\theta}$ such that $f_1(z) > \dots > f_D(z)$ then $0 < \underline{\delta}(G_1 > \dots > G_D)$. Also, if there exists $z \in H_{\theta}$ such that $f_1(z) \geq \dots \geq f_D(z)$ does not hold, then $\overline{\delta}(G_1 > \dots > G_D) > 1$. In particular, if there exists $y \in \R$ such that $F_1(y) > \dots > F_D(y)$, then $0 < \underline{\delta}(G_1 > \dots > G_D)$, and if there exists $n \in \Z$ such that $F_1(n) \geq \dots \geq F_D(n)$ does not hold, then $\overline{\delta}(G_1 > \dots > G_D) < 1$.
\end{prop}

\begin{rmq}
In the infinite-dimensional case, the only known lower bounds on the density are shown in particular cases by using delicate combinatorial arguments (\textit{cf.} \cite[2.2]{RuSa} and \cite[Theorem 2.5.1 (1)]{Dev}).
\end{rmq}

\subsection{Moments in the discrete case}

\label{Mom1}

To bound the probabilities involved in our results in the case $D=2$, moment estimates, such as Chebyshev's inequality, can prove to be very useful. In particular, races where, for each $a$, $\mathbb P(f_1(\nu_{\theta}^a Z_{\theta}) > f_2(\nu_{\theta}^a Z_{\theta})) \underset{q \to +\infty}{\longrightarrow} 1$ or $0$ ("extremely biased races") or each $\mathbb P(f(\nu_{\theta}^a Z_{\theta}) > f_2(\nu_{\theta}^a Z_{\theta})) \underset{q \to +\infty}{\longrightarrow} \frac{1}{2}$ ("moderately biased races") can be obtained from sufficiently good estimates on the corresponding means and variances as in \cite{FiJo}. Here the limits are taken with respect to a parameter $q$ attached to the prime number races considered. To go beyond the first two moments, and for example have an explicitly computable characteristic function at our disposal as in \cite[Theorem 3.4]{Cha} or \cite[Theorem 2.4]{ChaIm}), we would need to assume extra linear independence.\\

In the next Proposition, we give formulas for the first two moments of the random variables $f_j(\nu_{\theta}^a Z_{\theta})$ for some particular types of functions $f_j$.

\begin{prop} \label{Moments} Let $f = c + \sum_{k=1}^r a_k X_k + \overline{a_k} X_k^{-1} \in \C(X_1, \dots, X_r)$. Let $\theta_1, \dots, \theta_r$ be real numbers such that $\theta_i \not \in \pi\Q$ for $i \leq n$, and $\theta_i \in \pi\Q$ for $n < i \leq r$. Then for $0 \leq a \leq d-1$, one has $$\E(f(\nu_{\theta}^a Z_{\theta}))=c + \sum_{n < k \leq r} (a_k e^{ia \theta_k} + \overline{a_k} e^{-i a \theta_k})$$ and $$\Var(f(\nu_{\theta}^a Z_{\theta})) = 2\sum_{1 \leq k \leq n} |a_k|^2 + 4 \Re\left(\sum_{\underset{\theta_i + \theta_j \in \pi\Q}{1 \leq i < j \leq n}} a_i a_j e^{ia(\theta_i + \theta_j)}\right) + 4 \Re\left(\sum_{\underset{\theta_i - \theta_j \in \pi\Q}{1 \leq i < j \leq n}} a_i \overline{a_j} e^{ia(\theta_i-\theta_j)}\right).$$ In particular if $\theta_i \not \in \pi\Q$ for $1 \leq i \leq r$, then $\E(f(\nu_{\theta}^a Z_{\theta}))=c$ does not depend on $a$, and if no relation of the form $\theta_i \pm \theta_j \in \pi\Q$ holds for $1 \leq i < j \leq n$, then $\Var(f(\nu_{\theta}^a Z_{\theta})$ does not depend on $a$.
\end{prop}

\begin{demo} We note that the $n$ first components of $Z_{\theta}$ are products of non-zero integral powers of independent uniform random variables on $\T$, so they are uniform on $\T$, while the last $r-n$ components of $Z_{\theta}$ are (almost surely) equal to $1$. The result for the mean is then straightforward. For the variance, we expand the squared modulus and remark that $Z_{\theta, i} Z_{\theta, j} = 1$ (almost surely) if and only if $\theta_i + \theta_j \in \pi \Q$, and $Z_{\theta, i} \overline{Z_{\theta, j}} = 1$ (almost surely) if and only if $\theta - \theta_j \in \pi \Q$.
\end{demo}

As a corollary of the above Proposition, we deduce another sufficient condition for ties between those functions to have density zero.

\begin{cor} For $1 \leq j \leq D$, let $f_j = c_j + \sum_{k=1}^r a_k^{(j)} X_k + \overline{a_k^{(j)}} X_k^{-1} \in \C(X_1, \dots, X_r)$. Let $\theta_1, \dots, \theta_r$ be real numbers such that $\theta_i \not \in \pi\Q$ for $i \leq n$, and $\theta_i \in \pi\Q$ for $n < i \leq r$. Finally, let $G_j : t \mapsto f_j\left(e^{i\theta_1 t}, \dots, e^{i \theta_r t}\right) + o(1)$ as $t \to +\infty$.

\begin{itemize}
\item[i)] If $n > 1$ (i.e. we are in the non-degenerate case), and if  for every $a \in \{0, \dots, d-1\}$, the map $$j \mapsto \sum_{1 \leq k \leq n} |a_k^{(j)}|^2 + 2 \Re\left(\sum_{\underset{\theta_i + \theta_k \in \pi\Q}{1 \leq i < k \leq n}} a_i^{(j)} a_k^{(j)} e^{ia(\theta_i + \theta_k)}\right) + 2 \Re\left(\sum_{\underset{\theta_i - \theta_k \in \pi\Q}{1 \leq i < k \leq n}} a_i^{(j)} \overline{a_k^{(j)}} e^{ia(\theta_i-\theta_k)}\right)$$ is injective then for every permutation $\sigma$ of $\{1, \dots, D\}$, $\delta(G_{\sigma(1)} > \dots > G_{\sigma(D)})$ exists.
\item[ii)] If for every $a \in \{0, \dots, d-1\}$, the map $j \mapsto c_j + \sum_{k=n}^r \left(a_k^{(j)} e^{ia \theta_k} + \overline{a_k^{(j)}} e^{-i a \theta_k}\right) $ is injective then for every permutation $\sigma$ of $\{1, \dots, D\}$, $\delta(G_{\sigma(1)} > \dots > G_{\sigma(D)})$ exists.
\end{itemize}
\end{cor}

\begin{demo} 
\begin{itemize}
\item[i)] The hypothesis implies that $\Var(f_i(\nu_{\theta}^aZ_{\theta})) \neq \Var(f_j(\nu_{\theta}^aZ_{\theta}))$ for any $a \in \{0, \dots, d-1\}$ and any distinct $i, j \in \{1, \dots, D\}$ by Proposition \ref{Moments}, so that the random variable $f_i(\nu_{\theta}^aZ_{\theta}) - f_j(\nu_{\theta}^aZ_{\theta})$ is almost surely non-zero. By Lemma \ref{Atom}, this implies $\P(f_i(\nu_{\theta}^aZ_{\theta}) = f_j(\nu_{\theta}^aZ_{\theta}))=0$, and the result follows from Theorem \ref{Reste} \textit{ii)}.

\item[ii)] The proof is similar, except that this time $\E(f_i(\nu_{\theta}^aZ_{\theta}) - f_j(\nu_{\theta}^aZ_{\theta})) \neq 0$, which implies again that $\P(f_i(\nu_{\theta}^aZ_{\theta}) = f_j(\nu_{\theta}^aZ_{\theta}))=0$.
\end{itemize}
\end{demo}

\subsection{Moments in the infinite-dimensional case}

\label{Momente}

We work in the setting of Section \ref{Infini}. As before, we would like to have at least the first two moments of the random variable $S$ at our disposal. This is possible at the cost of assuming additional hypotheses on the coefficients $a_n$ of the function $F$, and on the linear relations between the almost periods of $F$, \textit{i.e.} the elements of the support $\Lambda(F)$ of $F$. We give below two distinct hypotheses that allow us to compute the second moment.\\

Let $\theta = (\theta_n)_{n \geq 1}$ be a sequence of pairwise distinct positive real numbers and for any $N \geq 1$ let $\Theta_N = \{\theta_n \mid n \leq N\}$ and $\Theta = \bigcup_{N \geq 1} \Theta_N$. We are going to do an analysis close to what we did in the previous sections, but we alter our notations to take into account the infinite number of $\theta_n$'s. Define inductively $\mathcal{B}_1 := \{\theta_1\}$ and for $N \geq 2$, $$\mathcal{B}_N := \left\{\begin{matrix}
\mathcal{B}_{N-1} \text{ if } \theta_N \in \Span_{\Q}(\Theta_{N-1})\\ \mathcal{B}_{N-1} \cup \{\theta_N\} \text{ otherwise.}
\end{matrix}\right.$$ We let $\mathcal{B} := \bigcup_{N \geq 1} \mathcal{B}_N$, so that $\mathcal{B}$ is a basis of $\Span_{\Q} \Theta$. For any $j \geq 1$, write the decomposition of $\theta_j$ in $\mathcal{B}$ as $$\theta_j = \sum_{\theta \in \mathcal{B}} c_{\theta, j} \theta$$ and let $d_j$ be the least common multiple of the denominators of the $c_{\theta, k}$ (written in irreducible form), for $\theta \in \mathcal{B}_j$ and $k \leq j$. Since the sequence of sets $(\mathcal{B}_j)_{j \geq 1}$ is increasing, we see that $d_j$ really only depends on $j$.

\begin{defi} Let $(Z_{\theta})_{\theta \in \mathcal B}$ be a sequence of independent random variables, uniform on $\T$. For any $N \geq 1$ and $n \leq N$, let $$Z_{n, N} = \prod_{\theta \in \mathcal{B}} Z_{\theta}^{d_N c_{\theta, n}}.$$
\end{defi}

Notice that, by definition, $d_N c_{n, \theta} \in \Z$ for $n \leq N$, and that for every $n \geq 1$, $c_{\theta, n} = 0$ for all but finitely many $\theta \in \mathcal{B}$, so the above product is a finite product.\\




\begin{lem} Let $c \in \C$, $(a_n)_{n \geq 1} \in \ell^2(\C)$ and for every $N \geq 1$, $$S_{\Theta_N} := c + \sum_{n \leq N} \left(a_n Z_{n, N} + \overline{a_n} \overline{Z_{n, N}}\right).$$ Assume that no $\theta_n$ is an integer multiple of another, that is for every $i, j \geq 1$ with $i \neq j$, we have $\theta_i \not \in \theta_j \Z$. Then $(S_{\Theta_N})_{N \geq 1}$ converges in $L^2$.
\end{lem}

\begin{demo} Since $L^2$ is complete, it is enough to prove that $(S_{\Theta_N})_{N \geq 1}$ is Cauchy in $L^2$. Let $m > n \geq 1$, then $$\E(|S_{\Theta_m} - S_{\Theta_n}|^2) = \E\left(\left|\sum_{k=1}^m (a_kZ_{k, m} + \overline{a_k} \overline{Z_{k, m}}) - \sum_{k=1}^n (a_kZ_{k, n} + \overline{a_k} \overline{Z_{k, n}}) \right|^2\right).$$

Expanding the square, we end up with terms of twelve different kinds : $a_ka_j\E(Z_{k,p}Z_{j,p})$, $a_k \overline{a_j}\E(Z_{k, p} \overline{Z_{j, p}})$ for $p \in \{m, n\}$ and $1 \leq j, k \leq p$, $-a_ka_j\E(Z_{k,m}Z_{j,n}), -a_k \overline{a_j}\E(Z_{k, m} \overline{Z_{j, n}})$ for $1 \leq k \leq m, 1 \leq j \leq n$ and their conjugates. Each $Z_{k, p}$ is uniform on $\T$ and the product of two uniform random variables on $\T$ is either $1$ or uniform on $\T$, in which case it has mean zero. Thus, it is enough to detect in which of the above cases we end up with $1$.

\begin{itemize}
\item[•] For $p \in \{m, n\}$ and $1 \leq j,k \leq p$ we have $$Z_{k,p} Z_{j,p} = \prod_{\mathcal \theta \in \mathcal{B}} Z_{\theta}^{d_p(c_{\theta, k} + c_{\theta, j})}.$$ Since the $Z_{\theta}, \theta \in \mathcal{B}$ are independent, this product is (almost surely) $1$ if and only if $d_p(c_{\theta, k} + c_{\theta, j})=0$ for every $\theta \in \mathcal{B}$. By definition, this means that $\theta_j = -\theta_k$ which can't be because each $\theta_n$ is positive. So those terms contribute $0$.

\item[•] Similarly, we have $$Z_{k,p} \overline{Z_{j,p}} = \prod_{\mathcal \theta \in \mathcal{B}} Z_{\theta}^{d_p(c_{\theta, k} - c_{\theta, j})}.$$ This is equal to $1$ when $j=k$, but when $j \neq k$ there is at least one non-zero exponent since the $\theta_n$'s are pairwise distinct.

\item[•] Now for $1 \leq k \leq m, 1 \leq j \leq n$, $$Z_{k,m}Z_{j,n} = \prod_{\theta \in \mathcal{B}} Z_{\theta}^{d_m c_{\theta, k} + d_n c_{\theta, j}}.$$ As before, this is $1$ if and only if $\theta_k = -\frac{d_n}{d_m} \theta_j$, which can't be since the $\theta_n$'s are positive.

\item[•] Finally for $1 \leq k \leq m, 1 \leq j \leq n$, $$Z_{k,m} \overline{Z_{j,n}} = \prod_{\theta \in \mathcal{B}} Z_{\theta}^{d_m c_{\theta, k} - d_n c_{\theta, j}},$$ and this is $1$ if and only if $\theta_j = \frac{d_m}{d_n} \theta_k$. But by definition, $d_n$ divides $d_m$ for $m > n$, so by hypothesis on the $\theta_n$'s the previous equality can only happen if and only if $k=j$.
\end{itemize}

Gathering everything, we have $$\E(|S_{\Theta_m} - S_{\Theta_n}|^2) = 2 \sum_{k=1}^m |a_k|^2 + 2 \sum_{k=1}^n |a_k|^2 - 4 \sum_{k=1}^n |a_k|^2 = 2 \sum_{k=n+1}^m |a_k|^2.$$

Since $(a_n)_{n \geq 1} \in \ell^2(\C)$, this proves that $(S_{\Theta_N})_{N \geq 1}$ is Cauchy for the $L^2$ norm.
\end{demo}

If the sequence $(d_j)_{j \geq 1}$ is bounded, it is stationary since it is non-decreasing. In that case, we let $d$ be its limit and $N_0 \geq 1$ be such that $d_N = d$ for every $N \geq N_0$. Then $Z_{n, N} = Z_{n, N_0}$ for any $N \geq \max(N_0, n)$, so for any $n \geq 1$ we let $Z_n := Z_{n, N_0}$. With these notations, we can now state the following result.

\begin{lem} Let $c \in \C$, $(a_n)_{n \geq 1} \in \ell^2(\C)$ and for every $N \geq 1$, $$S_{\Theta_N} := c + \sum_{n \leq N} \left(a_n Z_{n, N} + \overline{a_n} \overline{Z_{n, N}}\right).$$ Assume that $(d_j)_{j \geq 1}$ is bounded. Then $(S_{\Theta_N})_{N \geq 1}$ converges in $L^2$.
\end{lem}

\begin{demo} The proof is made easier in this case by the fact that $S_{\Theta_m} - S_{\Theta_n} = \sum_{k=n+1}^m (a_kZ_k + \overline{a_k} \overline{Z_k})$, and the terms of the sum are easily seen to be pairwise orthogonal. We conclude exactly as in the previous proof.
\end{demo}

\begin{cor} \label{L2} Let $F \sim c + \sum_{n \geq 1} \left(a_n e^{i \theta_n t} + \overline{a_n} e^{-i \theta_n t}\right)$ be $B^1$-almost periodic (see Theorem \ref{Helly}), where $c \in \C$ and $(a_n)_{n \geq 1} \in \ell^2(\C)$. Assume either that no $\theta_i$ is an integer multiple of another or that $(d_j)_{j \geq 1}$ is bounded. Then we can choose $$S_{\Theta} := \lim_{N \to +\infty} c + \sum_{n \leq N} \left(a_n Z_{n, N} + \overline{a_n} \overline{Z_{n, N}}\right)$$ (the limit being taken in $L^2$ norm) in the conclusion of Theorem \ref{Helly}, \textit{i.e.} for every bounded continuous function $g$ on $\R$, one has $$\frac{1}{X} \int_0^X g(F(y)) \,\mathrm{d}y \underset{X \to +\infty}{\longrightarrow} \E(g(S_{\Theta})).$$ Moreover, we have $\E(S_{\Theta}) = c$ and $\Var(S_{\Theta}) = 2 \sum_{n \geq 1} |a_n|^2.$
\end{cor}

\begin{demo} The first part is an immediate consequence of Theorem \ref{Helly} and the previous two lemmas. The last two formulas are straightforward.
\end{demo}

\begin{rmq}\begin{itemize}
\item[i)] The formula for the variance is classical when $F$ is $B^2$-almost periodic (\cite[p.109]{Bes}), but a $B^1$-almost periodic function with $\ell^2$ coefficients is not necessarily $B^2$ (though it is in the same $B^1$-class as at least one $B^2$-almost periodic function by \cite[p.110]{Bes}). 

\item[ii)] It is tempting to say that $(S_{\Theta_N})_{N \geq 1}$ converges to $S_{\Theta}$ in $L^1$ under no other hypothesis than $F$ being $B^1$-almost periodic, since $\limsup_{X \to +\infty} \frac{1}{X} \int_0^X |F(y) - P_N(y)| \,\mathrm{d}y$ goes to zero as $N$ tends to infinity. However we don't know if the previous quantity equals $\E(|S_{\Theta} - S_{\Theta_N}|)$.

\item[iii)] If the $\theta_n$'s are linearly independent over $\Q$, so that $\mathcal{B} = \Theta$, then the $Z_n$'s are pairwise independent, and one can prove the almost sure convergence of $\left(S_{\Theta_N}\right)_{N \geq 1}$ by Kolmogorov's two series theorem for example (see \cite[Theorem 2.5.6]{Dur}). In general, the almost sure convergence does not immediatly follow from the assumption that the series for $F(t)$ converges for every $t \in \R$ : the set $\Gamma := \left\{\left(e^{i \theta t}\right)_{\theta \in \mathcal B} \mid t \in \R\right\}$ is easily seen to be dense in $\T^{\mathcal B}$ as a consequence of the continuous version of the Kronecker-Weyl theorem, but it has measure zero. By the Riesz-Fischer theorem, we at least know that $\left(S_{\Theta_N}\right)_{N \geq 1}$ admits an almost surely converging subsequence.
\end{itemize}
\end{rmq}

\subsection{Prime divisor races over global function fields}

\label{2}

We now give an application of the previous results to the study of prime divisor races over global function field extensions.\\

Let $L/K$ be a geometric Galois extension of function fields, with constant field $\F_q$, the finite field with $q$ elements, \textit{i.e.} $K$ is a finitely generated extension of $\F_q$, has transcendence degree $1$ over $\mathbb F_q$, $\mathbb F_q$ is algebraically closed in $L$ (and therefore in $K$) and $L/K$ is Galois. We let $g_K$ and $g_L$ denote the genus (\cite[Theorem 6.6]{Ros}) of $K$ and $L$ respectively. Let $C_1, \dots, C_D$ be $D \geq 1$ distinct conjugacy classes of $G := \Gal(L/K)$. Define $$\pi_{C_i}(n) := \#\{P \text{ prime divisor of } K \text{ unramified in } L \mid \deg(P)=n, \Frob_P=C_i\},$$ where $\Frob_P$ denotes the Frobenius conjugacy class of $P$ in $G$. The Chebotarev density theorem (\cite[Theorem 9.13B]{Ros}) states that $$\pi_{C_i}(n) = \frac{|C_i|}{|G|} \frac{q^n}{n} + O\left(\frac{q^{n/2}}{n}\right).$$ This shows that $$\frac{\pi_{C_i}(n)}{|C_i|} \underset{n \to +\infty}{\sim} \frac{\pi_{C_j}(n)}{|C_j|}$$ for any two $i \neq j$, but we want to compare those two quantities beyond this first order asymptotic. The question is, how often can it happen that $$\frac{\pi_{C_1}(n)}{|C_1|} >\dots > \frac{\pi_{C_D}(n)}{|C_D|} \,?$$ More precisely, we are interested in the following density, provided it exists, $$\delta(L/K; C_1, \dots, C_D) := \lim_{X \to +\infty} \frac{\#\left\{n \leq X \ \vline\  \frac{\pi_{C_1}(n)}{|C_1|} > \dots >  \frac{\pi_{C_D}(n)}{|C_1|}\right\}}{X}.$$ 


When studying the densities $\delta(L/K; C_{\sigma(1)}, \dots, C_{\sigma(D)})$ for every permutation $\sigma$ of $\{1, \dots, D\}$, we say we study the prime divisor race between $C_1, \dots, C_D$. As usual, we will denote by $\overline{\underline{\delta}}(L/K; C_1, \dots, C_D)$ the corresponding $\liminf$ and $\limsup$.

To study the above densities, we use the Artin $L$-functions associated to irreducible characters of $G$. If $\chi$ is such a character, one has the following convenient expression : \begin{equation}\label{Artin}\log L(s, \chi) = \sum_P \sum_{n \geq 1} \frac{\chi(P^n) q^{-n \deg(P) s}}{n}\end{equation} for $\Re(s) > 1$, where $\chi(P^n)$ is a short way of writing $\chi(\Frob_P^n)$ when $P$ is unramified, and it is $$\frac{1}{e(P)} \sum_{g \in I(\mathfrak{P})} \chi(g (\Frob_P)^n)$$ when $P$ is ramified. Here, $\Frob_P$ is any preimage in $G$ of the corresponding residual Frobenius automorphism, $e(P)$ is the ramification index of $P$ in $L$ and $I(\mathfrak{P})$ is the inertia subgroup of a prime $\mathfrak{P}$ of $L$ dividing $P$. Since $I(\mathfrak{P})$ and $\Frob_P$ only depends on $P$ up to conjugacy and $\chi$ is a central function on $G$, this is indeed independent of any choice.

It is convenient to move from the variable $s$ to the variable $u := q^{-s}$, and to write $\mathcal{L}(u, \chi) := L(s, \chi)$. Then one has the following important theorem (\cite[Theorems 9.16A and 9.16B]{Ros}) :

\begin{thm}[Weil] \label{Func} The function $L(s, \chi_0)$, where $\chi_0$ is the trivial character of $G$ (which is also the zeta function $\zeta_K$ of $K$) is a rational function in $u$ with integer coefficients, which we factorize as \begin{equation}\label{zeta}\zeta_K(s) = \mathcal L(u, \chi_0) = \frac{\prod_{j=1}^{2g_K} (1-\gamma(\chi_0, j) u)}{(1-u)(1-qu)}.\end{equation} If $\chi \neq \chi_0$ is a non-trivial irreducible character of $G$, then $\mathcal{L}(u, \chi)$ is a polynomial in $u$ with integer coefficients which we factorize as \begin{equation}\label{L}\mathcal{L}(u, \chi) = \prod_{j=1}^{M_{\chi}} (1-\gamma(\chi, j)u)\end{equation} for some integer $M_{\chi} \geq 0$. The $\gamma(\chi ,j)$ are called the inverse zeros of $\mathcal{L}(u, \chi)$ and have absolute value $\sqrt{q}$ (Riemann Hypothesis for curves over $\F_q$). Moreover, if $\gamma$ is an inverse zero of $\mathcal{L}(u, \chi)$ then $\frac{q}{\gamma} = \overline{\gamma}$ is an inverse zero of $\mathcal{L}(u, \overline{\chi})$.
\end{thm}

The last statement of the theorem is a simple consequence of the functional equation satisfied by Artin $L$-functions, which we do not formulate here (for non-trivial characters, one has to combine the functional equation of Hecke $L$-series and Brauer's induction theorem as in the case of Artin $L$-functions over number fields \cite[p.81-83]{Langl}).\\

We now make some preliminary work to study the prime divisor races in $L/K$. We let $C$ be a conjugacy class of $G$.

On the one hand, by the definition (\ref{Artin}) of Artin $L$-functions, one has \begin{equation} \label{f1}u \frac{\mathrm{d}}{\mathrm{d}u} \log \mathcal{L}(u, \chi) = \sum_P \sum_{n \geq 1} \deg P \chi(P^n) u^{n \deg P} = \sum_{n \geq 1} \left(\sum_{\underset{\deg P \mid n}{P}} \deg P \chi(P^{\frac{n}{\deg P}})\right) u^n.\end{equation}

On the other hand, from the above factorizations (\ref{zeta}) and (\ref{L}), we obtain \begin{equation} \label{f2} u \frac{\mathrm{d}}{\mathrm{d}u} \log \mathcal{L}(u, \chi_0) = \sum_{n \geq 1} \left(q^n + 1 - \sum_{j=1}^{2g_K} \gamma(\chi, j)^n\right) u^n\end{equation} and for $\chi \neq \chi_0$ \begin{equation} \label{f3}u \frac{\mathrm{d}}{\mathrm{d}u} \log \mathcal{L}(u, \chi) = \sum_{n \geq 1} \left(- \sum_{j=1}^{M_{\chi}} \gamma(\chi, j)^n\right) u^n.\end{equation}

Writing $$u \frac{\mathrm{d}}{\mathrm{d}u} \log \mathcal{L}(u, \chi) = \sum_{n \geq 1} c_n(\chi) u^n,$$ and using the second orthogonality relations on characters, we find from formula (\ref{f1}) that \begin{align*}\sum_{\chi} \overline{\chi(C)} c_n(\chi) &= \sum_{\underset{\deg P \mid n}{P}} \deg P \sum_{\chi} \chi(P^{\frac{n}{\deg P}}) \overline{\chi(C)}\\ &= \frac{\#G}{\#C} \sum_{d \mid n} d \#\{P \mid \deg P = d, \Frob_P^{\frac{n}{d}} \subset C\}\\ &= \frac{\#G}{\#C} n \pi_C(n) + R_C(n) + O\left(q^{n/3}\right),
\end{align*} where $$R_C(n) = \left\{\begin{matrix}
\frac{n}{2} \frac{\#G}{\#C} \#\{P \mid \deg P = \frac{n}{2}, \Frob_P^2 \subset C\} \text{ if } n \text{ is even}\\ 0 \text{ otherwise}.
\end{matrix}\right.$$

Using now formulas (\ref{f2}) and (\ref{f3}), we obtain $$\sum_{\chi} \overline{\chi(C)} c_n(\chi) = q^n + 1 - \sum_{j=1}^{2g_K} \gamma(\chi_0, j)^n - \sum_{\chi \neq \chi_0} \overline{\chi(C)} \sum_{j=1}^{M_{\chi}} \gamma(\chi, j)^n.$$

Combining the two expressions for $\sum_{\chi} \overline{\chi(C)} c_n(\chi)$ we obtain \begin{equation}\label{pi}\frac{\#G}{\#C} \pi_C(n) = \frac{q^n}{n} - \frac{R_C(n)}{n} - \frac{1}{n} \sum_{j=1}^{2g_K} \gamma(\chi_0, j)^n - \frac{1}{n} \sum_{\chi \neq \chi_0} \overline{\chi(C)} \sum_{j=1}^{M_{\chi}} \gamma(\chi, j)^n + O\left(\frac{q^{n/3}}{n}\right).\end{equation} We now introduce $$C^{1/2} := \{g \in G \mid g^2 \in C\},$$ and remark it is stable by conjugacy in $G$, so it is the disjoint union of conjugacy classes $D_1, \dots, D_t$ of $G$. Moreover, \begin{align*}\#\{P \mid \deg P = \frac{n}{2}, \Frob_P^2 \subset C\} = \sum_{i=1}^t \pi_{D_i}\left(\frac{n}{2}\right) &= \sum_{i=1}^t \frac{\#D_i}{\#G} \frac{2}{n} q^{n/2} + O(q^{n/4})\\ &= \frac{\#(C^{1/2})}{\#G} \frac{2}{n} q^{n/2} + O(q^{n/4})\end{align*} by the above formula.

This shows that \begin{align*}R_C(n) &= \left\{\begin{matrix}
\frac{\#(C^{1/2})}{\#C} q^{n/2} + O(q^{n/4}) \text{ if } n \text{ is even}\\ 0 \text{ otherwise}
\end{matrix}\right.\\ &= \frac{\#(C^{1/2})}{2 \#C} q^{n/2} + \frac{\#(C^{1/2})}{2 \#C} q^{n/2} e^{i \pi n} + O(q^{n/4}).\end{align*}

Finally, replacing $R_C(n)$ in (\ref{pi}), we obtain \begin{align*}\frac{\#G}{\#C} \pi_C(n) &= \frac{q^n}{n} - \frac{\#(C^{1/2})}{\#C} \frac{q^{n/2}}{2n} - \frac{\#(C^{1/2})}{\#C} \frac{q^{n/2}}{2n} e^{i \pi n} - \frac{1}{n} \sum_{j=1}^{2g_K} \gamma(\chi_0, j)^n\\ & - \frac{1}{n} \sum_{\chi \neq \chi_0} \overline{\chi(C)} \sum_{j=1}^{M_{\chi}} \gamma(\chi, j)^n + O\left(\frac{q^{n/3}}{n}\right).\end{align*}

Similarly, with $\pi_K(n) := \#\{P \mid \deg P = n\}$, one has $$\pi_K(n) = \frac{q^n}{n} - \frac{q^{n/2}}{2n} - e^{i \pi n} \frac{q^{n/2}}{2n} - \frac{1}{n} \sum_{j=1}^{2g_K} \gamma(\chi_0, j)^n + O\left(q^{n/3}\right).$$

Combining those two formulas, we obtain $$\frac{n}{q^{n/2}} \left(\frac{\#G}{\#C} \pi_C(n) - \pi_K(n)\right) = \frac{1 - \frac{\#(C^{1/2})}{\#C}}{2}  + \frac{1 - \frac{\#(C^{1/2})}{\#C}}{2} e^{i \pi n} - \sum_{\chi \neq \chi_0} \overline{\chi(C)} \sum_{j=1}^{M_{\chi}} \left(\frac{\gamma(\chi, j)}{\sqrt q}\right)^n + o(1)$$ as $n \to +\infty$.


Grouping pairs of conjugate inverse zeros we have shown :

\begin{prop} \label{Expli} Let $\gamma_1, \dots, \gamma_r$ be the inverse zeros with positive imaginary part of the $\mathcal L(u, \chi)$, for $\chi \neq \chi_0$, counted without multiplicity. For $1 \leq j \leq r,$ write $\gamma_j = \sqrt{q} e^{i\theta_j}$. Then for any conjugacy class $C$ of $G$ we have $$\frac{n}{q^{n/2}} \left(\frac{\#G}{\#C} \pi_C(n) - \pi_K(n)\right) = r_C + z_C + a_{\pi}(C) e^{i \pi n} - \sum_{j=1}^r \left(a_j(C) e^{i \theta_j n} + \overline{a_j(C)} e^{- i \theta_j n}\right) + o(1)$$ as $n \to +\infty$, where $$r_C := \frac{1 - \frac{\#(C^{1/2})}{\#C}}{2},$$ $$z_C := - \sum_{\chi \neq \chi_0} \overline{\chi(C)} \ord_{u = q^{-1/2}} \mathcal{L}(u, \chi),$$ $$a_{\pi}(C) = r_C - \sum_{\chi \neq \chi_0} \overline{\chi(C)}\ord_{u = -q^{-1/2}} \mathcal{L}(u, \chi),$$ and for $1 \leq j \leq r$, $$a_j(C) := \sum_{\chi \neq \chi_0} \overline{\chi(C)} \ord_{u = \gamma_j^{-1}} \mathcal{L}(u, \chi).$$
\end{prop}

We have thus shown that a suitable rescaling of $\pi_C(n)$ is of the form we studied in the previous sections. The rescaling of $\frac{\pi_C(n)}{\#C}$ does not depend on $C$, so that will allow us to study prime divisor races between conjugacy classes of $G$.

\begin{thm} \label{MainF} Let $C_1, \dots, C_D$ be conjugacy classes of $G$. For $1 \leq j \leq D$, let $f_j = r_{C_j} + z_{C_j} + a_{\pi}(C_j) \frac{X_{r+1} + X_{r+1}^{-1}}{2} - \sum_{k=1}^r \left(a_k(C_j) X_j + \overline{a_k(C_j)} X_j^{-1}\right)\in \C(X_1, \dots, X_{r+1})$ and $F_j : t \mapsto f_j\left(e^{i \theta_1 t}, \dots, e^{i \theta_r t}, e^{i \pi t}\right)$.

\begin{itemize}
\item[i)] Degenerate case : Assume $\theta_i \in \pi\Q$ for $1 \leq i \leq r$, i.e. that each $\mathcal{L}(u, \chi)$, $\chi \neq \chi_0$, is a product of (rescaled) cyclotomic polynomials. If there exists $n \in \Z$ such that $F_1(n) > \dots > F_D(n)$, then $0 < \underline{\delta}(L/K; C_1, \dots, C_D)$ and if there exists $n \in \Z$ such that $F_1(n) \geq \dots \geq F_D(n)$ does not hold, then $\overline{\delta}(L/K; C_1, \dots, C_D) < 1$. Moreover, if for $1 \leq j \leq D-1$, and for $0 \leq n \leq d-1$, one has $F_j(n) \neq F_{j+1}(n)$, then $\delta(L/K; C_1, \dots, C_D)$ exists.

\item[ii)] Non-degenerate case : Assume $\theta_i \not \in \pi \Q$ for at least one $i \in \{1, \dots, r\}$. If there exist $a \in \{0, \dots, d-1\}$ and $z \in H_{\theta}$ such that $f_1(\nu_{\theta}^az) > \dots > f_D(\nu_{\theta}^az)$ then $0 < \underline{\delta}(L/K; C_1, \dots, C_D)$. Also, if there exist $a \in \{0, \dots, d-1\}$ and $z \in H_{\theta}$ such that $f_1(\nu_{\theta}^az) \geq \dots \geq f_D(\nu_{\theta}^az)$ does not hold, then $\overline{\delta}(L/K; C_1, \dots, C_D) < 1$. In particular, if there exists $n \in \Z$ such that $F_1(n) > \dots > F_D(n)$, then $0 < \underline{\delta}(L/K; C_1, \dots, C_D)$, and if there exists $n \in \Z$ such that $F_1(n) \geq \dots \geq F_D(n)$ does not hold, then we have $\overline{\delta}(L/K; C_1, \dots, C_D) < 1$. Moreover, if for $0 \leq a \leq d-1$ and $1 \leq j \leq D-1$, there exists $n \equiv a \emph{ mod } d$ such that $F_j(n) \neq F_{j+1}(n)$, then $\delta(L/K; C_1, \dots, C_D)$ exists.
\end{itemize}
\end{thm}

\begin{demo} This is an immediate application of Proposition \ref{Expli}, Theorem \ref{Reste} and Proposition \ref{PosDisc}.
\end{demo}

We now treat an example for which there is linear dependence between the $\theta_i$'s. This example was featured in \cite{ChaIm}, but since it did not satisfy the required linear independence condition under which the authors worked, they weren't able to study the corresponding prime divisor race. Take $K = \F_7(t)$ and $L = K(\alpha)$ where $\alpha$ has minimal polynomial $f = X^6 - (t^2+t)X^3 - 1$ over $K$. Then as detailed in \cite[4.2]{ChaIm}, $G = \Gal(L/K) \simeq \mathfrak{S}_3$. We note $C_1 = \{\text{id}\}$, $C_2 = \{(1 \, 2), (1 \, 3), (2 \, 3)\}$ and $C_3 = \{(1 \, 2 \, 3), (1 \, 3 \, 2)\}$, and it is well-known that we have the following character table for $G$ :

$$\begin{array}{c|c|c|c}
\mathfrak{S}_3 & C_1 & C_2 & C_3 \\ \hline
\chi_0 &1&1&1\\ \hline
\chi_1 & 1 & -1 & 1\\ \hline
\chi_2 & 2 &0 & -1
\end{array}$$

One has $$\mathcal{L}(u, \chi_1) = 1 + 4u + 7u^2 = (1-\gamma_1 u)(1 - \overline{\gamma_1}u),$$ $$\mathcal{L}(u, \chi_2) = 1 + u + 7u^2 = (1-\gamma_2 u)(1 - \overline{\gamma_2}u),$$ with (those two values are inverted in \cite{ChaIm}) $$\gamma_1 = -2 + i \sqrt 3$$ and $$\gamma_2 = \frac{-1 + 3i \sqrt 3}{2}.$$ Then we have $\theta_1 = \arctan\left(-\frac{\sqrt{3}}{2}\right)$, $\theta_2 = \arctan(-3\sqrt 3)$ and $\theta_1 + \theta_2 = \frac{4\pi}{3}$. Adding $\theta_3 = \pi$ because of the coefficient $a_{\pi}(C)$, we are in the non-degenerate case (because $\theta_2 = - \arccos \left(\frac{1}{\sqrt{28}}\right)$ as one easily verifies, and such a number is known not to be a rational multiple of $\pi$, see \cite{Var} for example). With notations from Section \ref{Disc}, we have $m=1, d=6, c_2 = \frac{2}{3}, b_{1, 2} = -1, c_3 = \frac{1}{2}, b_{1, 3} = 0$.

A quick computation using PARI/GP shows that for $i \neq j \in \{1,2,3\}$ and every $a \in \{0,1,2,3,4,5\}$, there exists $n \equiv a \text{ mod } 6$ such that $F_i(n) \neq F_j(n)$, so for every permutation $\sigma \in \mathfrak{S}_3$, the density $\delta(L/K; C_{\sigma(1)}, C_{\sigma(2)}, C_{\sigma(3)})$ exists, \textit{i.e.} the race between $C_1, C_2$ and $C_3$ is weakly inclusive (Definition \ref{Incl}). Also, for every permutation $\sigma \in \mathfrak{S}_3$, the inequality $F_{\sigma(1)}(n) > F_{\sigma(2)}(n) > F_{\sigma(3)}(n)$ happens for some $n \leq 7$ so we may conclude that $0 < \delta(L/K; C_{\sigma(1)}, C_{\sigma(2)}, C_{\sigma(3)}) < 1$ and in particular the race between $C_1, C_2$ and $C_3$ is inclusive (Definition \ref{Incl}).\\

\begin{rmq} If one wants to study races between functions counting prime divisors of degree less than $n$, instead of equal to $n$ as above, one can use the following explicit formula (\cite[Theorem 2.1]{ChaIm}) : \begin{align*}\frac{n}{q^{n/2}} \left(\frac{\#G}{\#C} \sum_{k=1}^n \pi_C(k) - \pi_K(n)\right) &= r_C \frac{q + \sqrt{q}}{q-1} + r_C \frac{q-\sqrt{q}}{q-1} e^{in\pi} - 2 \sum_{j=1}^{2g_K} \frac{\gamma(\chi_0,j)}{\gamma(\chi_0, j)-1} e^{i n \theta(\chi_0, j)}\\ &- \sum_{\chi \neq \chi_0} \overline{\chi(C)} \sum_{j=1}^{M_{\chi}} \frac{\gamma(\chi,j)}{\gamma(\chi, j)-1} e^{i n \theta(\chi, j)} + o(1)\end{align*} as $n \to +\infty$, and use our method similarly since this has the shape we studied above.
\end{rmq}

\section*{Appendix : a proof of the discrete Kronecker-Weyl theorem}

\begin{demoKW} Recall $\theta_1, \dots, \theta_r$ are real numbers such that $\{\pi, \theta_1, \dots, \theta_r\}$ is linearly independent over $\Q$, $$\Gamma = \left\{\left(e^{i\theta_1X}, \dots, e^{i\theta_rX}\right) \mid X \in \Z\right\}$$ and we want to show that for every continuous $f : \T^r \rightarrow \C$, $$\frac{1}{X} \sum_{n=1}^X f\left(e^{i\theta_1n}, \dots, e^{i\theta_rn}\right) \underset{X \to +\infty}{\longrightarrow} \int_{\T^r} f \,\mathrm{d}\lambda,$$ where $\lambda$ is the Lebesgue measure on $\T^r$.

By the Stone-Weierstrass theorem, it is enough to prove the result when $f$ is a trigonometric polynomial, that is, a linear combination of monomials in $z_1, \dots, z_r$. Indeed, if the result is true for such functions, then for any $\varepsilon > 0$, we can find a trigonometric polynomial $g$ such that $||f-g||_{\infty} < \varepsilon$, and for every $X$ large enough we have $$\left|\frac{1}{X} \sum_{n=1}^X g\left(e^{i\theta_1n}, \dots, e^{i\theta_rn}\right) - \int_{\T^r} g \,\mathrm{d}\lambda\right| < \varepsilon.$$ For such $X$, we find \begin{align*}\left|\frac{1}{X} \sum_{n=1}^X f\left(e^{i\theta_1n}, \dots, e^{i\theta_rn}\right) - \int_{\T^r} f \,\mathrm{d}\lambda\right| &\leq \left|\frac{1}{X} \sum_{n=1}^X (f-g)\left(e^{i\theta_1n}, \dots, e^{i\theta_rn}\right) - \int_{\T^r} (f-g) \,\mathrm{d}\lambda\right|\\ &+ \left|\frac{1}{X} \sum_{n=1}^X g\left(e^{i\theta_1n}, \dots, e^{i\theta_rn}\right) - \int_{\T^r} g \,\mathrm{d}\lambda\right| < 3\varepsilon\end{align*} which proves that $$\frac{1}{X} \sum_{n=1}^X f\left(e^{i\theta_1n}, \dots, e^{i\theta_rn}\right) \underset{X \to +\infty}{\longrightarrow} \int_{\T^r} f \,\mathrm{d}\lambda.$$ By linearity, we now only have to prove the theorem for monomials $$\begin{array}{rccl} &\T^r &\longrightarrow &\C\\f : &(z_1,\dots, z_r)&\mapsto & z_1^{n_1} \dots z_r^{n_r},\end{array}$$ where $n_1, \dots, n_r \in \Z$.

The result is obviously true if $(n_1, \dots, n_r) = (0, \dots, 0)$, \textit{i.e.} if $f=1$ since both sides are equal to $1$. Now assume at least one $n_i$ is non-zero. On the one hand we have $$\int_{\T^r} f \,\mathrm{d}\lambda = \prod_{k=1}^r \int_{\T} z^{n_k} \,\mathrm{d}\lambda = 0.$$ On the other hand, since $\{\pi, \theta_1, \dots, \theta_r\}$ is linearly independant over $\Q$, we have that $n_1\theta_1 + \dots + n_r\theta_r \not \in 2\pi\Z$, so that $e^{i(n_1\theta_1 + \dots + n_r \theta_r)} \neq 1$. Now, summing the geometric progression, we find \begin{align*}
\frac{1}{X} \sum_{n=1}^X f\left(e^{i\theta_1n}, \dots, e^{i\theta_rn}\right) &= \frac{1}{X} \sum_{n=1}^X e^{in(n_1\theta_1 + \dots + n_r \theta_r)}\\ &= \frac{1}{X} \frac{e^{i(X+1)(n_1\theta_1 + \dots + n_r \theta_r)} - e^{i(n_1\theta_1 + \dots + n_r \theta_r)}}{e^{i(n_1\theta_1 + \dots + n_r \theta_r)} - 1}\\ &\underset{X \to +\infty}{\longrightarrow} 0
\end{align*} since $\frac{e^{i(X+1)(n_1\theta_1 + \dots + n_r \theta_r)} - e^{i(n_1\theta_1 + \dots + n_r \theta_r)}}{e^{i(n_1\theta_1 + \dots + n_r \theta_r)} - 1}$ is bounded.
\end{demoKW}

\begin{rmq} The continuous version of the Kronecker-Weyl theorem states that, assuming $\theta_1, \dots, \theta_r$ are linearly independent over $\Q$, for every continuous function $f : \T^r \rightarrow \C$, one has $$\frac{1}{X} \int_0^X f\left(e^{i\theta_1y}, \dots, e^{i\theta_ry}\right) \,\mathrm{d}y \underset{X \to +\infty}{\longrightarrow} \int_{\T^r} f \,\mathrm{d}\mu.$$ Its proof is similar as the one given above. The first step reduces to the case of trigonometric polynomials by using the Stone-Weierstrass theorem, and the last calculation is done with integrals instead of discrete sums.
\end{rmq}

\subsection*{Acknowledgements}

The author would like to thank Lucile Devin, Daniel Fiorilli, Florent Jouve, Nathan Ng, Günter Rote and Nicolae Strungaru for useful discussions, and Peter Humphries for providing proofs of both versions of the Kronecker-Weyl theorem, the continuous one in a MathOverflow thread, and the discrete one in his master thesis \cite{Hum}.

Data sharing not applicable to this article as no datasets were generated or analysed during the current study. The author declares he has no conflict of interest.

\bibliographystyle{plain}
\nocite{*}
\bibliography{bibli}
\label{sec:name}
\addcontentsline{toc}{chapter}{References}

\bigskip
\footnotesize
\noindent Alexandre Bailleul, \textsc{ENS Paris-Saclay, Centre Borelli, UMR 9010, 91190 Gif-sur-Yvette, France}\\
\textit{Email address:} \textsf{alexandre.bailleul@ens-paris-saclay.fr}

\end{document}